\documentclass[reqno,centertags, 12pt]{amsart}
\usepackage{amsmath,amsthm,amscd,amssymb}
\usepackage{latexsym}
\usepackage{graphicx}

\newcommand{\bbA}{{\mathbb{A}}}
\newcommand{\bbC}{{\mathbb{C}}}
\newcommand{\bbD}{{\mathbb{D}}}
\newcommand{\bbE}{{\mathbb{E}}}
\newcommand{\bbP}{{\mathbb{P}}}
\newcommand{\bbR}{{\mathbb{R}}}
\newcommand{\bbZ}{{\mathbb{Z}}}

\newcommand{\calC}{{\mathcal C}}

\newcommand{\calL}{{\mathcal L}} 
\newcommand{\calM}{{\mathcal M}} 
\newcommand{\calV}{{\mathcal V}}

\newcommand{\lb}{\label}
\newcommand{\f}{\frac}
\newcommand{\ul}{\underline}
\newcommand{\ol}{\overline}
\newcommand{\ti}{\tilde  }
\newcommand{\wti}{\widetilde  }

\newcommand{\tr}{\text{\rm{Tr}}}

\newcommand{\s}{\text{\rm{s}}}

\newcommand{\supp}{\text{\rm{supp}}}

\newcommand{\bi}{\bibitem}

\newcommand{\beq}{\begin{equation}}
\newcommand{\eeq}{\end{equation}}
\newcommand{\ba}{\begin{align}}
\newcommand{\ea}{\end{align}}
\newcommand{\veps}{\varepsilon}

\DeclareMathOperator{\Real}{Re}

\DeclareMathOperator{\Ima}{Im}

\allowdisplaybreaks
\numberwithin{equation}{section}

\newtheorem{theorem}{Theorem}[section]

\newtheorem{proposition}[theorem]{Proposition}
\newtheorem{lemma}[theorem]{Lemma}

\theoremstyle{definition}

\newtheorem{example}[theorem]{Example}
\newtheorem{conjecture}[theorem]{Conjecture}

\theoremstyle{remark}

\newcommand{\abs}[1]{\lvert#1\rvert}

\newcounter{smalllist}
\newenvironment{SL}{\begin{list}{{\rm\roman{smalllist})}}{%
\setlength{\topsep}{0mm}\setlength{\parsep}{0mm}\setlength{\itemsep}{0mm}%
\setlength{\labelwidth}{2em}\setlength{\leftmargin}{2em}\usecounter{smalllist}%
}}{\end{list}} 

\newcommand{\bigtimes}{\mathop{\mathchoice%
{\smash{\vcenter{\hbox{\LARGE$\times$}}}\vphantom{\prod}}%
{\smash{\vcenter{\hbox{\Large$\times$}}}\vphantom{\prod}}%
{\times}%
{\times}%
}\displaylimits}

\begin{document}

\title[Fine Structure of the Zeros of OP, I]
{Fine Structure of the Zeros of Orthogonal Polynomials, \\I. A Tale of Two Pictures}
\author[B.~Simon]{Barry Simon*}

\thanks{$^*$ Mathematics 253-37, California Institute of Technology, Pasadena, CA 91125, USA. 
E-mail: bsimon@caltech.edu. Supported in part by NSF grant DMS-0140592} 
\thanks{Talk given at Constructive Functions Tech-04, Atlanta, November 2004}

\date{November 8, 2004}

\begin{abstract} Mhaskar-Saff found a kind of universal behavior for the 
bulk structure of the zeros of orthogonal polynomials for large $n$. Motivated 
by two plots, we look at the finer structure for the case of random Verblunsky 
coefficients and for what we call the BLS condition: $\alpha_n =Cb^n + 
O((b\Delta)^n)$. In the former case, we describe results of Stoiciu. In 
the latter case, we prove asymptotically equal spacing for the bulk of zeros. 
\end{abstract}

\maketitle

\section{Prologue: A Theorem of Mhaskar and Saff} \lb{s1} 

A recurrent theme of Ed Saff's work has been the study of zeros of orthogonal polynomials 
defined by measures in the complex plane. So I was happy that some thoughts I've had 
about zeros of orthogonal polynomials on the unit circle (OPUC) came to fruition  
just in time to present them as a birthday bouquet. To add to the appropriateness, the 
background for my questions was a theorem of Mhaskar-Saff \cite{MhS1} and the idea of 
drawing pictures of the zeros was something I learned from some of Ed's papers 
\cite{Saff,newLSS}. Moreover, ideas of Barrios-L\'opez-Saff \cite{BLS} played a role 
in the further analysis. 

Throughout, $d\mu$ will denote a probability measure on $\partial\bbD=\{z\in\bbC 
\mid \abs{z}=1\}$  which is nontrivial in that it is not supported on a finite set. 
$\Phi_n(z)$ (resp.~$\varphi_n(z)$) will denote the monic orthogonal polynomials 
(resp.~orthonormal polynomials $\varphi_n =\Phi_n/\|\Phi_n\|$). I will follow my 
book \cite{OPUC1,OPUC2} for notation and urge the reader to look there for 
further background. 

A measure is described by its Verblunsky coefficients 
\begin{equation} \lb{1.1a} 
\alpha_n = -\ol{\Phi_{n+1}(0)} 
\end{equation}
which enter in the Szeg\H{o} recursion 
\begin{align} 
\Phi_{n+1}(z) &= z\Phi_n (z) - \bar\alpha_n \Phi_n^*(z) \lb{1.1}  \\
\Phi_{n+1}^* (z) &= \Phi_n^* (z) - \alpha_n z\Phi_n(z)  \lb{1.2} 
\end{align}
where 
\begin{equation} \lb{1.3} 
P_n^*(z) = z^n \, \ol{P_n (1/\bar z)} 
\end{equation}

$\Phi_n$ has all its zeros in $\bbD$ \cite[Theorem~1.7.1]{OPUC1}. We let $d\nu_n$ be the 
pure point measure on $\bbD$ which gives weight $k/n$ to each zero of $\Phi_n$ of multiplicity 
$k$. For simplicity, we will suppose there is a $b\in [0,1]$ so that 
\begin{equation} \lb{1.4} 
\lim \abs{\alpha_n}^{1/n}=b 
\end{equation}
(root asymptotics). If $b=1$, we also need 
\begin{equation} \lb{1.5} 
\lim_{n\to\infty} \, \f{1}{n} \sum_{j=0}^{n-1} \, \abs{\alpha_j}=0 
\end{equation}
(which automatically holds if $b<1$). Here is the theorem of Mhaskar-Saff \cite{MhS1}: 

\begin{theorem}[Mhaskar-Saff \cite{MhS1}]\lb{T1.1} If \eqref{1.4} and \eqref{1.5} hold, 
then $d\nu_n$ converges weakly to the uniform measure on the circle of radius $b$. 
\end{theorem} 

We note that both this result and the one of Nevai-Totik \cite{NT89} I will mention in a moment 
define $b$ by $\limsup \abs{\alpha_n}^{1/n}$ and the Mhaskar-Saff result holds for $d\nu_{n(j)}$ 
where $n(j)$ is a subsequence with $\lim_{j\to\infty} \abs{\alpha_{n(j)}}^{1/n(j)}=b$. 

I want to say a bit about the proof of Theorem~\ref{T1.1} in part because I will need a slight 
refinement of the first part (which is from Nevai-Totik \cite{NT89}) and in part because I want to 
make propaganda for a new proof \cite{OPUC1} of the second part. 

The proof starts with ideas from Nevai-Totik \cite{NT89} that hold when $b<1$: 
\begin{SL} 
\item[(1)] By \eqref{1.1} and $\abs{\Phi_n (e^{i\theta})}=\abs{\Phi_n^* (e^{i\theta})}$, 
one sees inductively that 
\begin{equation} \lb{1.6} 
\sup_{n,z\in\partial\bbD}\, \abs{\Phi_n^*(z)} \leq \prod_{j=0}^\infty \, (1+\abs{\alpha_j})<\infty  
\end{equation}
and so, by the maximum principle and \eqref{1.3}, 
\begin{equation} \lb{1.7} 
C\equiv \sup_{n,z\notin\bbD}\, \abs{z}^{-n} \abs{\Phi_n(z)} <\infty
\end{equation}

\item[(2)] By \eqref{1.2}, 
\[
\sum_{n=0}^\infty \, \abs{\Phi_{n+1}^*(z) - \Phi_n^*(z)} \leq C\sum_{n=0}^\infty \, 
\abs{\alpha_n} \, \abs{z}^{n+1} <\infty 
\]
if $\abs{z}b<1$ by \eqref{1.4}. Thus, in the disk $\{z\mid \abs{z}<b^{-1}\}$, $\Phi_n^*$ has 
a limit. Since $b<1$, the Szeg\H{o} condition (see \cite[Chapter 2]{OPUC1}) holds, so 
\begin{equation} \lb{1.8} 
\varphi_n^*(z) \to D(z)^{-1} 
\end{equation}
on $\bbD$ (see \cite[Theorem~2.4.1]{OPUC1}), we conclude that $D(z)^{-1}$ has an analytic 
continuation to the disk of radius $b^{-1}$ and \eqref{1.8} holds there. (Nevai-Totik also 
prove a converse: If the Szeg\H{o} condition holds, $d\mu_\s =0$ and $D(z)^{-1}$ has an 
analytic continuation to the disk of radius $b^{-1}$, then $\limsup \abs{\alpha_n}^{1/n}  
\leq b$.) 

\item[(3)] When \eqref{1.4} holds, $D(z)$ is analytic in $\bbD$ and continuous on $\bar\bbD$ 
so $D(z)^{-1}$ has no zeros in $\bar\bbD$. We define the Nevai-Totik points $\{z_j\}_{j=1}^N$ 
($N$ in $\{0,1,2,\dots,\}\cup\infty$) with $1>\abs{z_1} \geq \abs{z_2} \geq \cdots > b$ to 
be all the solutions of $D(1/\bar z)^{-1}=0$ in $\bbA_b =\{z\mid b<\abs{z}<1\}$. Since \eqref{1.8} 
holds and $\varphi_n^*(z)=0\Leftrightarrow \varphi_n (1/\bar z)=0$, Hurwitz's theorem implies 
that the $z_j$ are precisely the limit points of zeros of $\varphi_n$ in the region $\bbA_b$. 
If $N =\infty$, $\abs{z_j}\to b$ so we conclude that for each $\veps >0$, the number of zeros of 
$\varphi_n$ in $\{z\mid\abs{z} > b+\veps\}$ is bounded as $n\to\infty$. 
\end{SL}

\smallskip
That concludes our summary of some of the results from Nevai-Totik. The next step is from 
Mhaskar-Saff. 

\begin{SL} 
\item[(4)] By \eqref{1.1a}, if $\{z_{jn}\}_{j=1}^n$ are the zeros of $\varphi_n(z)$ counting 
multiplicity, then 
\[
\abs{\alpha_n} = \prod_{j=1}^n \, \abs{z_{jn}} 
\]
so, by \eqref{1.4}, 
\[
\lim_{n\to\infty} \, \f{1}{n} \sum_{j=1}^n \log \abs{z_{jn}}=\log b 
\]
This together with (3) implies that the bulk of zeros must asymptotically lie 
on the circle of radius $b$ and, in particular, any limit point of $d\nu_n$ 
must be a measure on $\{z\mid \abs{z}=b\}$. 
\end{SL} 

\smallskip
Mhaskar-Saff complete the proof by using potential theory to analyze the 
limit points of the $\nu_n$. Instead, I will sketch a different idea from 
\cite[Section~8.2]{OPUC1} that exploits the CMV matrix (see \cite{CMV} and 
\cite[Sections~4.2--4.5]{OPUC1}). 

\begin{SL}
\item[(5)] If $\calC^{(N)}$ is the $N\times N$ matrix in the upper left of $\calC$, then 
\begin{equation} \lb{1.9a}
\Phi_N(z) = \det (z-\calC^{(N)}) 
\end{equation}
and, in particular, 
\begin{equation} \lb{1.9b}
\int z^k\, d\nu_N(z) =\f{1}{N}\, \tr ([\calC^{(N)}]^k) 
\end{equation}
It is not hard to see that \eqref{1.5} implies that on account of the structure 
of $\calC$, the right side of \eqref{1.9b} goes to $0$ as $N\to\infty$ for each 
$k>0$. Thus any limit point, $d\nu$,  of $d\nu_N$ has $\int z^k\, d\nu(z)=0$ for 
$k=1,2,\dots$. That determines the measure $d\nu$ uniquely since the only measure 
on $b[\partial\bbD]$ with zero moments is the uniform measure.  
\end{SL} 

\smallskip
That completes the sketch of the proof of the Mhaskar-Saff theorem. Before going 
on, I have two remarks to make. It is easy to see (\cite[Theorem~8.2.6]{OPUC1}) 
that $\langle \delta_\ell, \calC^k \delta_\ell \rangle = \int_0^{2\pi} e^{ik\theta} 
\abs{\varphi_\ell (e^{i\theta})}^2 \, d\mu(\theta)$. Thus the moments of the 
limit points of the Ce\`saro average 
\[
d\eta_N=\f{1}{N} \sum_{j=0}^{N-1}\, \abs{\varphi_j (e^{i\theta})}^2\, d\mu(\theta) 
\]
are the same as the moments of the limits of $d\nu_N$ (so if $d\nu_N$ has a limit 
that lives on $\partial\bbD$, $d\eta_N$ has the same limit). 

Second, a theorem like Mhaskar-Saff holds in many other situations. For example, 
if $\beta_n$ is periodic and $\alpha_n -\beta_n\to 0$, then the $d\nu_n$ converge 
to the equilibrium measure for the essential support of $d\mu$, which is a finite 
number of intervals. (See Ed Saff's book with Totik \cite{SaffTot}.) 

\smallskip 
One critical feature of the Mhaskar-Saff theorem is its universality. So long as we 
look at cases where \eqref{1.4} and \eqref{1.5} hold, the angular distribution is the 
same. Our main goal here is to go beyond the universal setup where the results will 
depend on more detailed assumptions on asymptotics. In particular, we will want to 
consider two stronger conditions than root asymptotics, \eqref{1.4}, namely, ratio 
asymptotics 
\begin{equation} \lb{1.8x} 
\lim_{n\to\infty} \, \f{\alpha_{n+1}}{\alpha_n}  =b
\end{equation} 
for some $b\in (0,1)$ and what we will call BLS asymptotics (or the BLS condition): 
\begin{equation} \lb{1.9} 
\alpha_n = Cb^n + O((b\Delta)^n)  
\end{equation}
where $b\in (0,1)$, $C\in\bbC$, and $\Delta\in (0,1)$. 

The name BLS is for Barrios-L\'opez-Saff \cite{BLS} who proved

\begin{theorem}[Barrios-L\'opez-Saff \cite{BLS}]\lb{T1.2} A set of Verblunsky coefficients 
obeys the BLS condition if and only if $d\mu_\s =0$ and the numeric inverse of the 
Szeg\H{o} function $D(z)^{-1}$, defined initially for $z\in\bbD$, has a meromorphic continuation to a disk of radius 
$(b\Delta')^{-1}$ for some $\Delta'\in (0,1)$ which is analytic except for a single pole 
at $z=b^{-1}$. 
\end{theorem}

This is Theorem~2 in \cite{BLS}. They allow $b\in\bbD\backslash\{0\}$ but, by the rotation 
invariance of OPUC (see \cite[Eqns.~(1.6.62)--(1.6.67)]{OPUC1}), any $b=\abs{b}e^{i\theta}$ 
can be reduced to $\abs{b}$. Another proof of Theorem~\ref{T1.2} can be found in 
\cite[Section 7.2]{OPUC1} where it is also proven that 
\begin{equation} \lb{1.10} 
C = -\bigl[\, \lim_{z\to b^{-1}}\, (z^{-1}-b)D(z)^{-1}\bigr] \, \ol{D(b)}
\end{equation}
and that if $\Delta$ or $\Delta'$ is in $(\f12, 1)$, then $\Delta =\Delta'$. 

We summarize the contents of the paper after the next and second preliminary 
section. I'd like to thank Mourad Ismail, Rowan Killip, Paul Nevai, and 
Mihai Stoiciu for useful conversations.

\section{Two Pictures and an Overview of the Fine Structure} \lb{s2} 

Take a look at two figures (Figure~1 and Figure~2) from my book \cite[Section~8.4]{OPUC1}. 
The first shows the zeros of $\Phi_n(z)$ for $\alpha_n = (\f12)^{n+1}$ and the second 
for $\alpha_n$ independent random variables uniformly distributed on $[-\f12,\f12]$. 
(Of course, the second is for a particular choice of the random variables made by 
Mathematica using its random number generator.) They are shown for $n=5,10,20,50,100$, 
and $200$. 

\begin{center}
\begin{figure}
\mbox{}\hfill\hbox to 0mm{$n=5$\hss}\hspace{50mm}\hfill\hspace{50mm}\hbox to 0mm{\hss$n=10$}\hfill\mbox{}\\[-2ex]
\mbox{}\hfill\includegraphics[scale=.50]{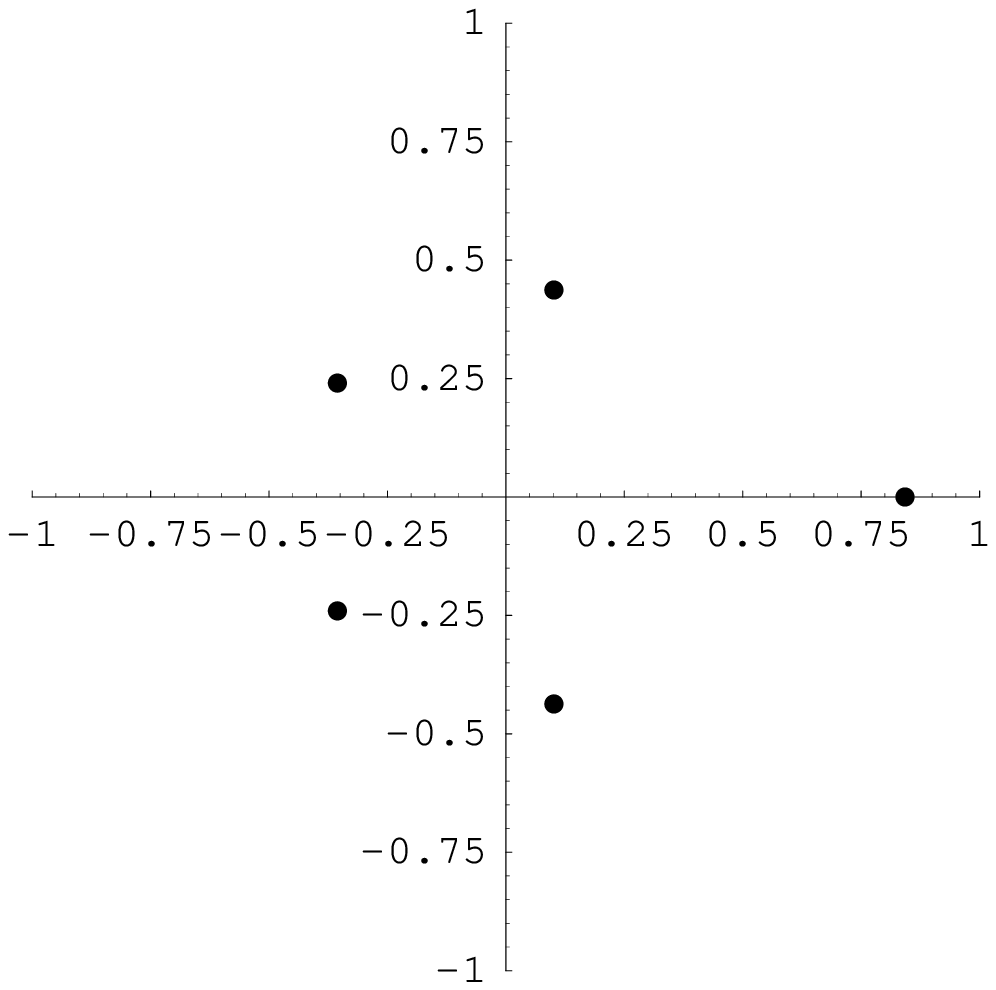}  \hfill\includegraphics[scale=.50]{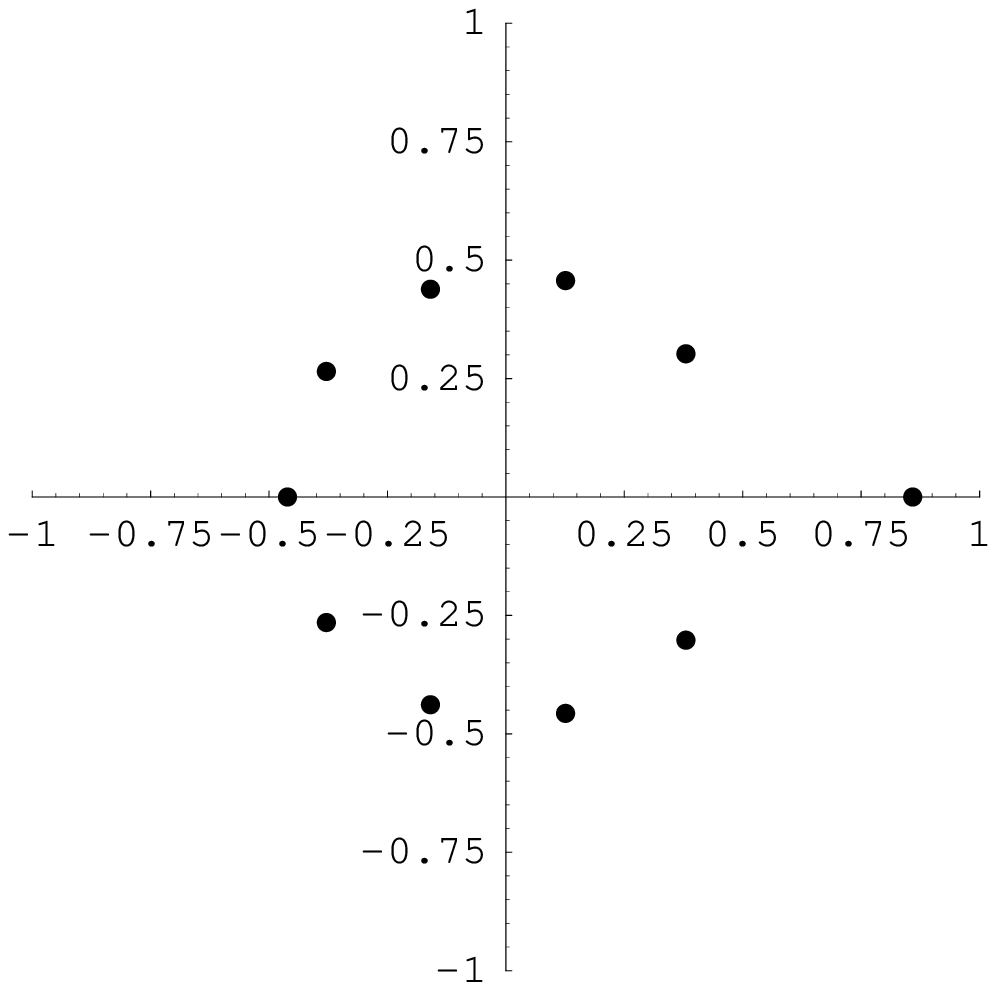}\hfill\mbox{}\\[8mm]
\mbox{}\hfill\hbox to 0mm{$n=20$\hss}\hspace{50mm}\hfill\hspace{50mm}\hbox to 0mm{\hss$n=50$}\hfill\mbox{}\\[-2ex]
\mbox{}\hfill\includegraphics[scale=.50]{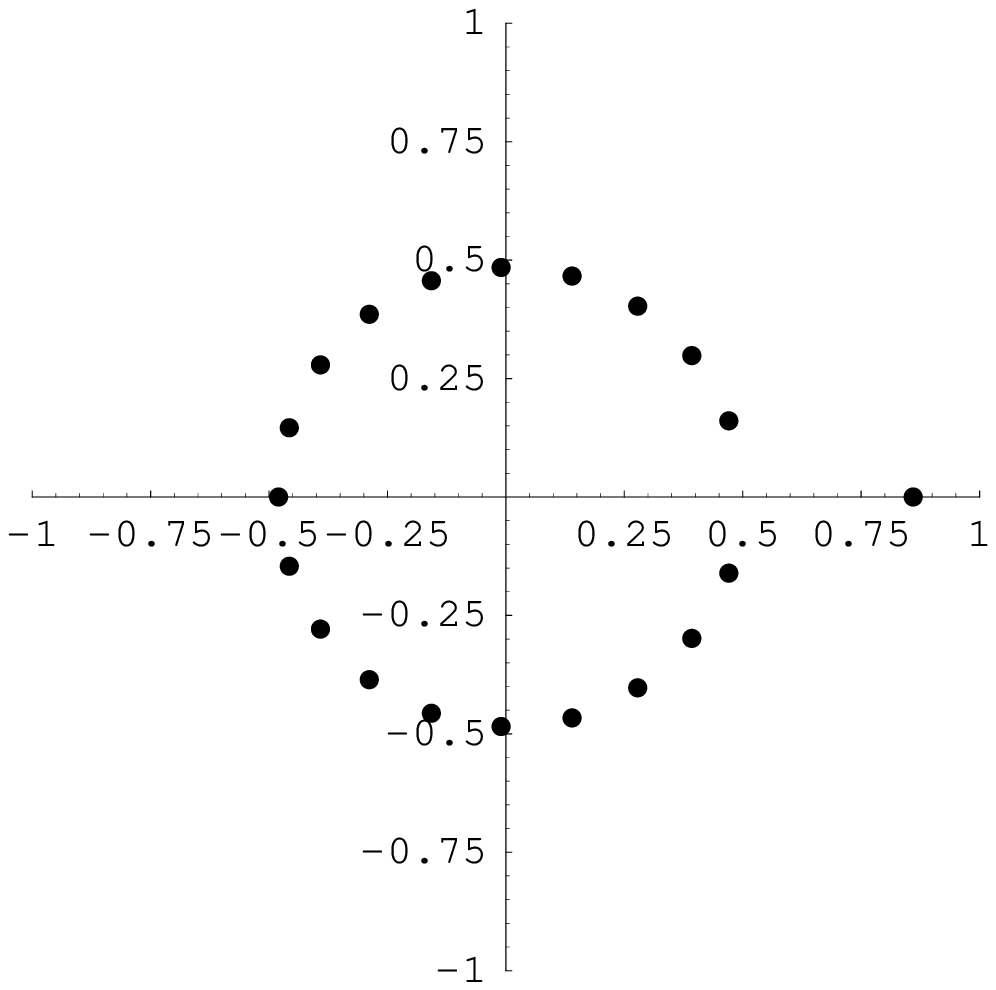} \hfill\includegraphics[scale=.50]{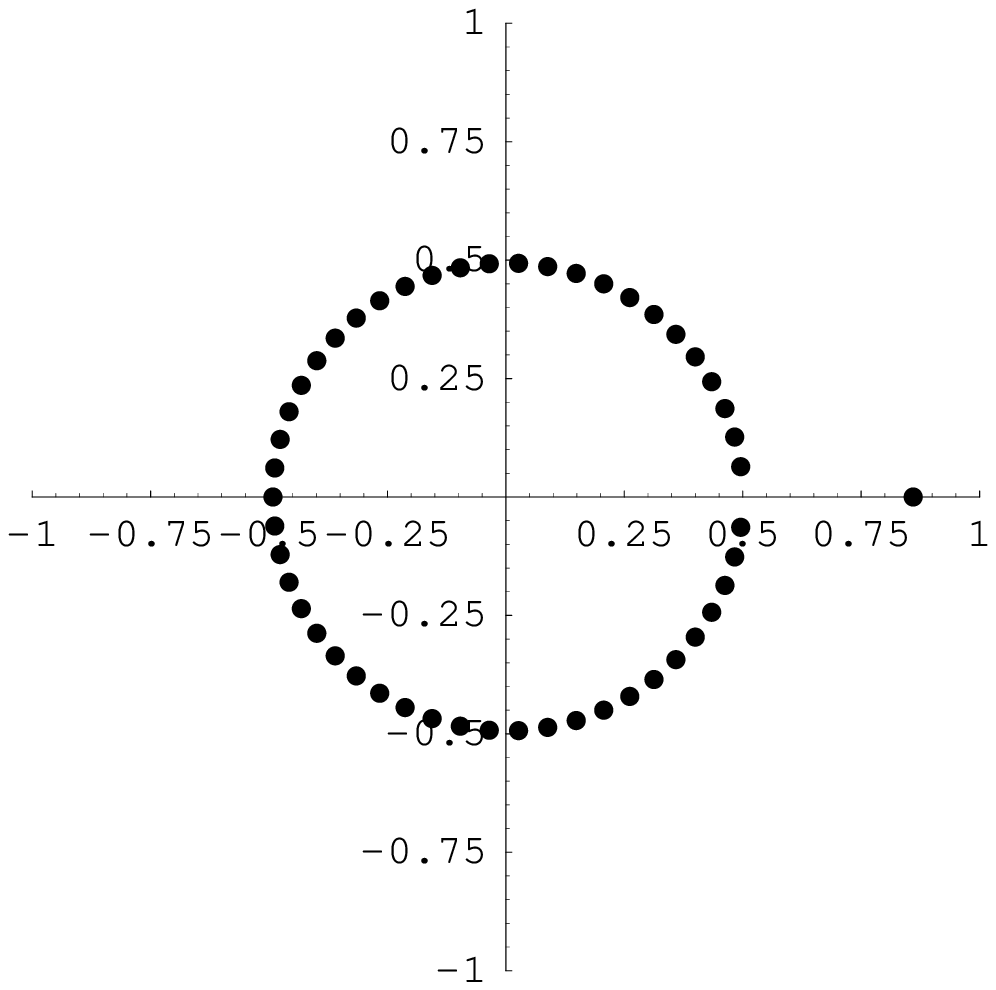}\hfill\mbox{}\\[8mm] 
\mbox{}\hfill\hbox to 0mm{$n=100$\hss}\hspace{50mm}\hfill\hspace{50mm}\hbox to 0mm{\hss$n=200$}\hfill\mbox{}\\[-2ex]
\mbox{}\hfill\includegraphics[scale=.50]{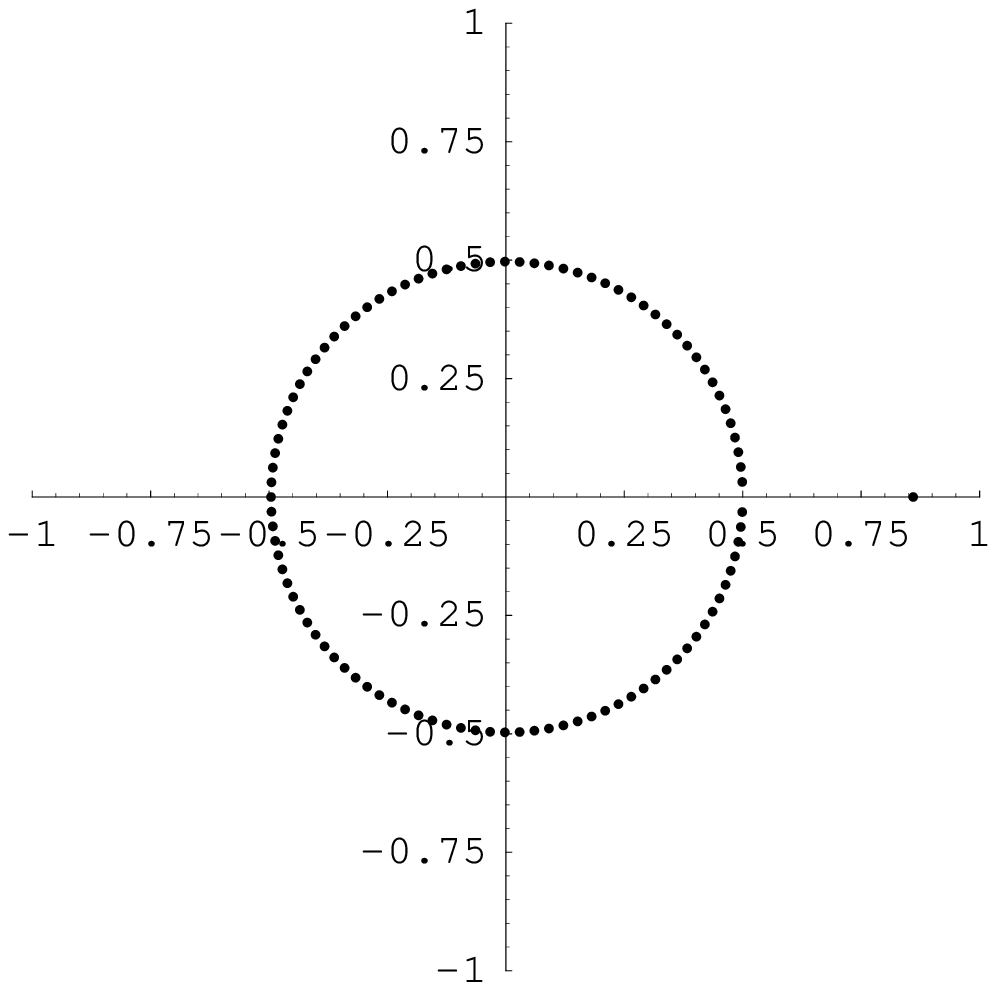}\hfill\includegraphics[scale=.50]{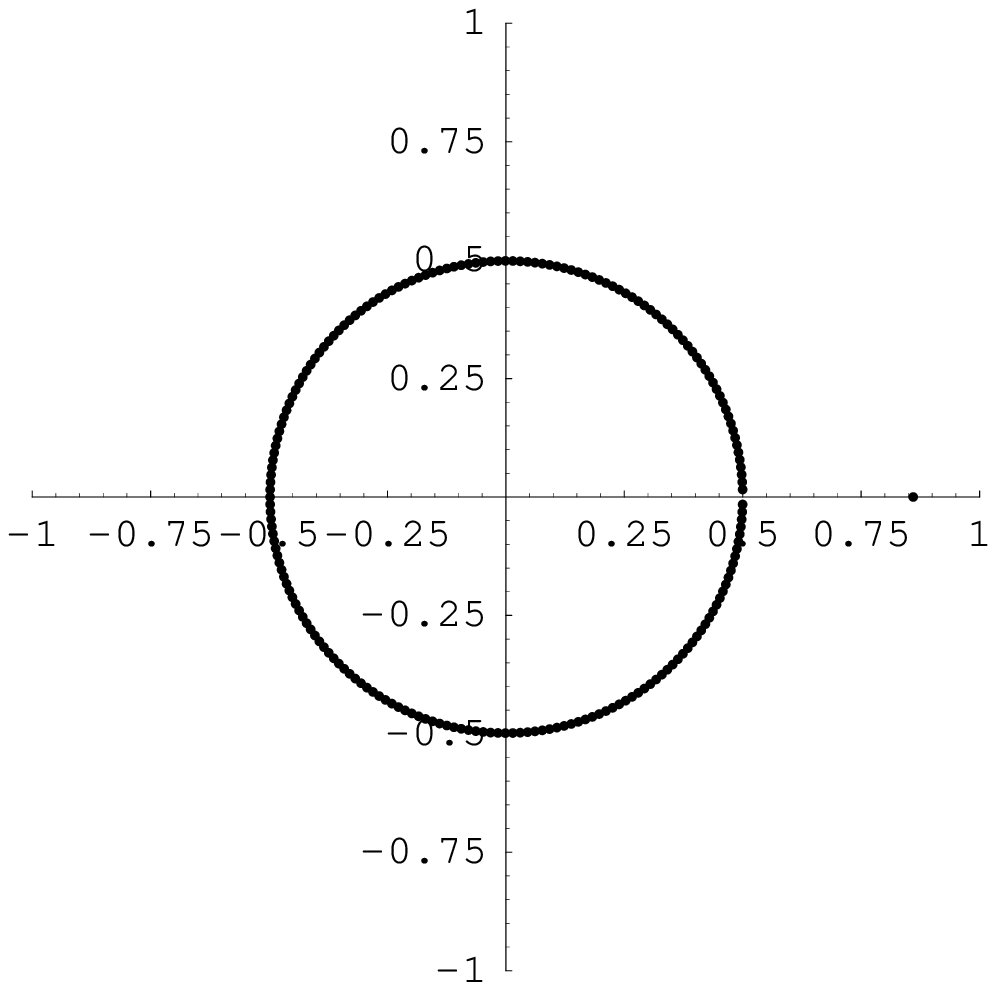}\hfill\mbox{}\\
\caption{}
\end{figure}
\end{center} 

\newpage

\begin{center}
\begin{figure}[t]
\mbox{}\hfill\hbox to 0mm{$n=5$\hss}\hspace{50mm}\hfill\hspace{50mm}\hbox to 0mm{\hss$n=10$}\hfill\mbox{}\\[-2ex]
\mbox{}\hfill\includegraphics[scale=.50]{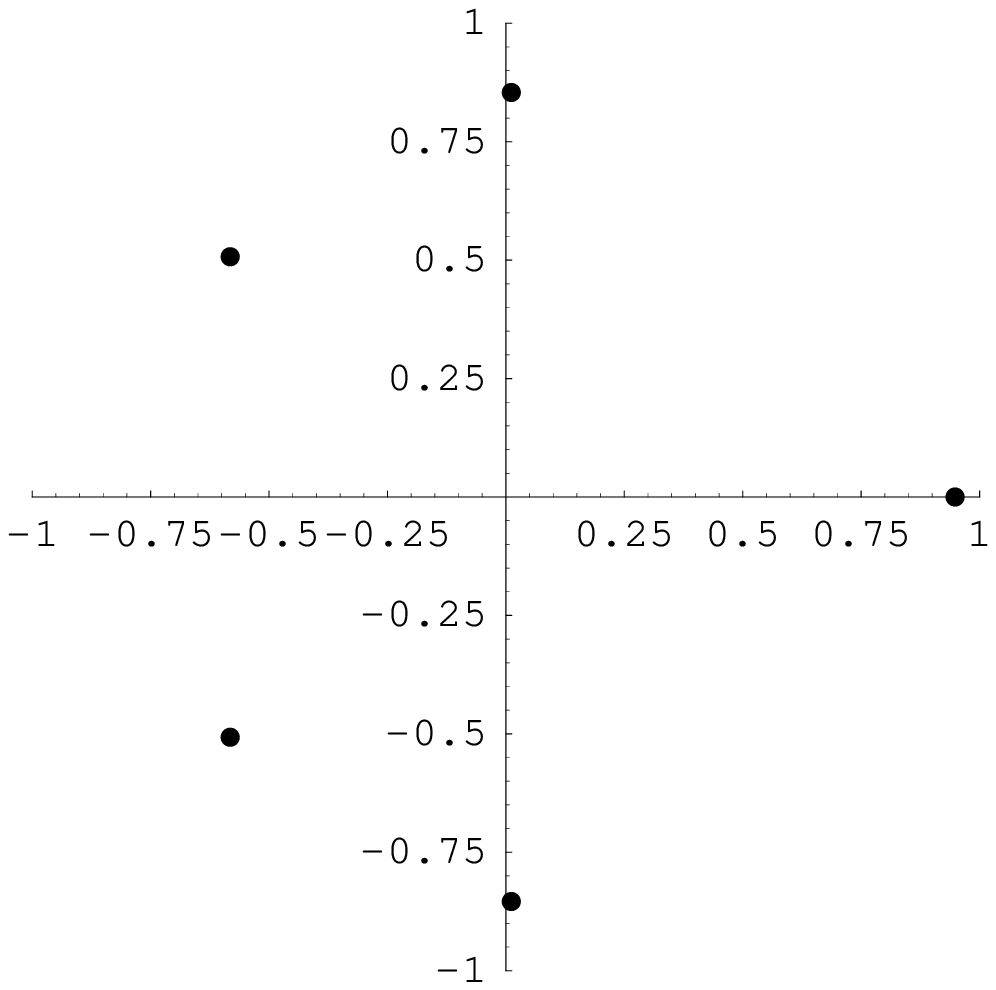}  \hfill\includegraphics[scale=.50]{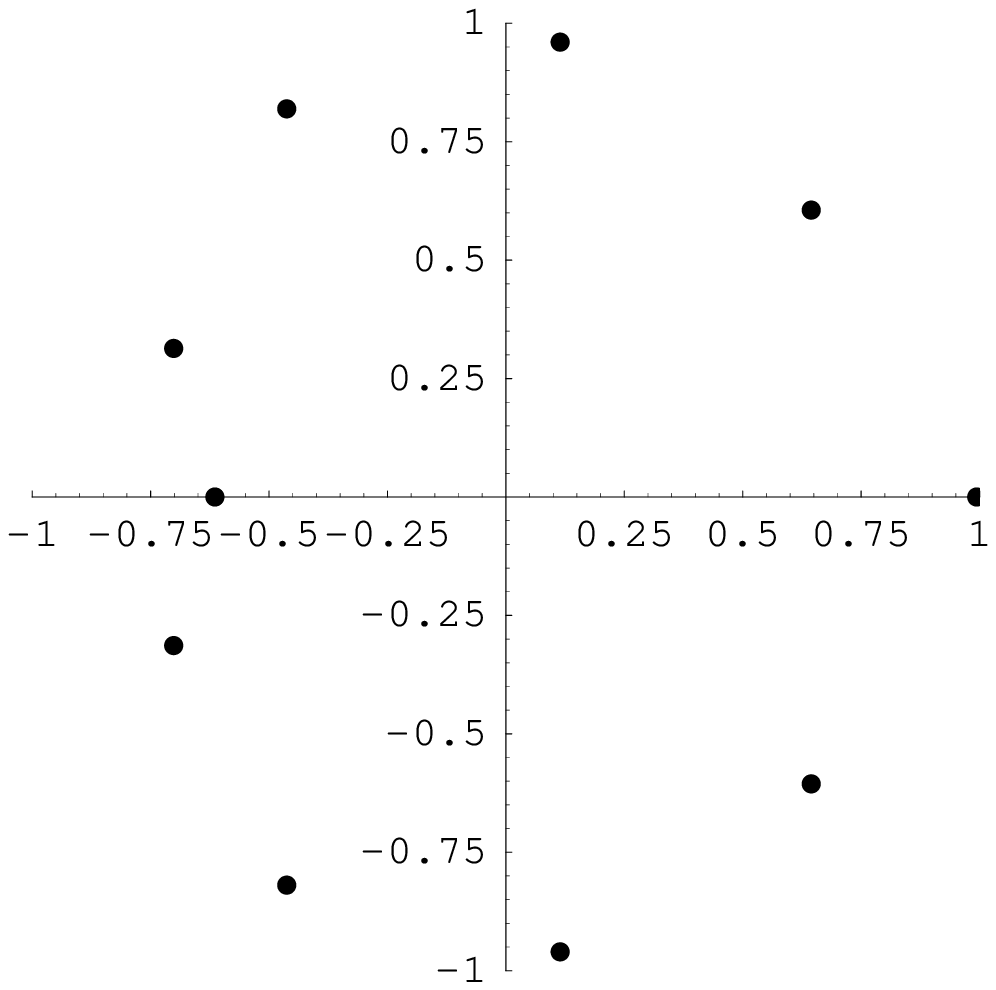}\hfill\mbox{}\\[8mm]
\mbox{}\hfill\hbox to 0mm{$n=20$\hss}\hspace{50mm}\hfill\hspace{50mm}\hbox to 0mm{\hss$n=50$}\hfill\mbox{}\\[-2ex]
\mbox{}\hfill\includegraphics[scale=.50]{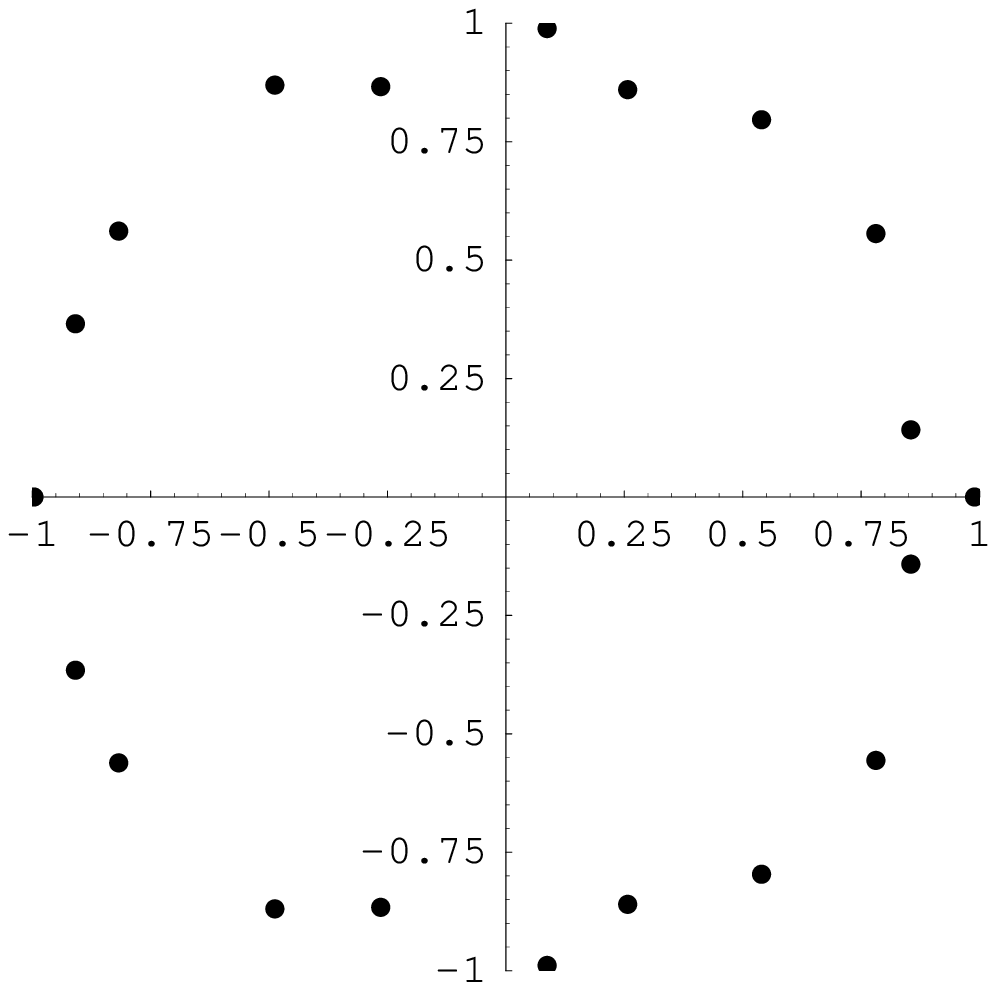} \hfill\includegraphics[scale=.50]{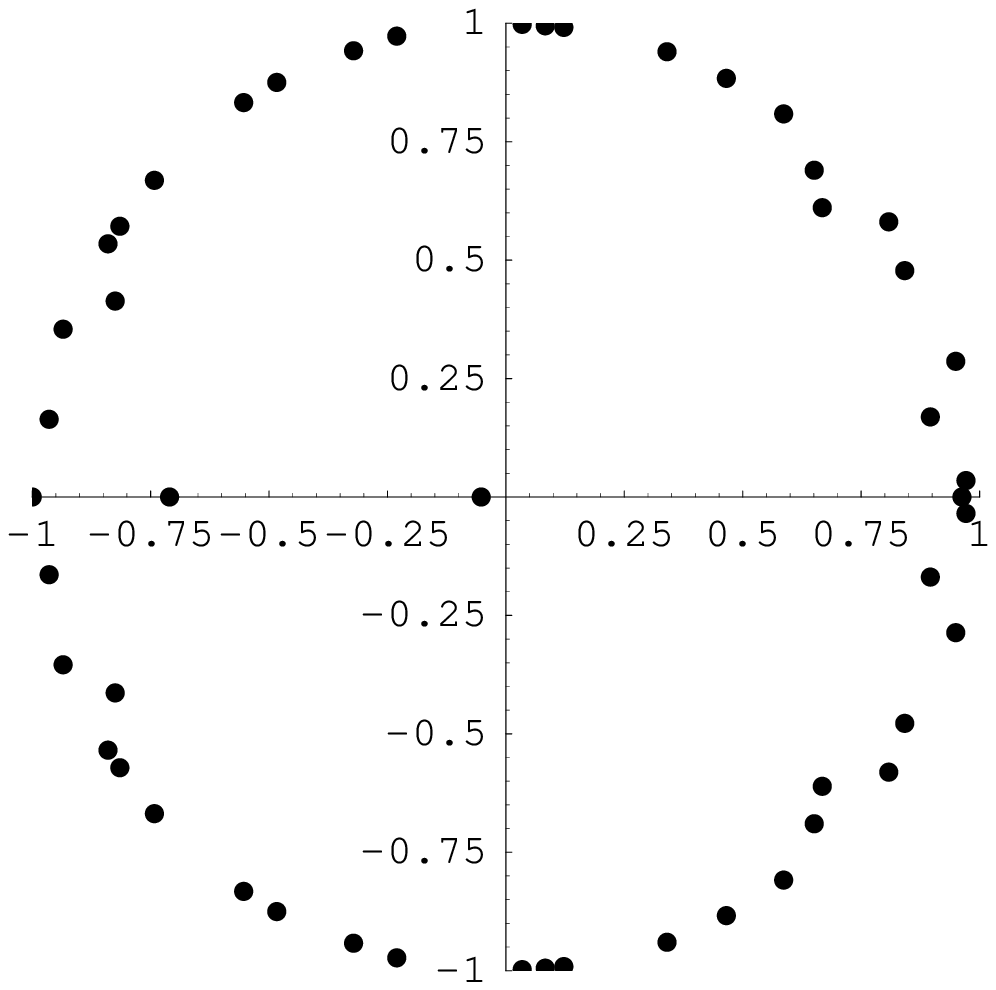}\hfill\mbox{}\\[8mm] 
\mbox{}\hfill\hbox to 0mm{$n=100$\hss}\hspace{50mm}\hfill\hspace{50mm}\hbox to 0mm{\hss$n=200$}\hfill\mbox{}\\[-2ex]
\mbox{}\hfill\includegraphics[scale=.50]{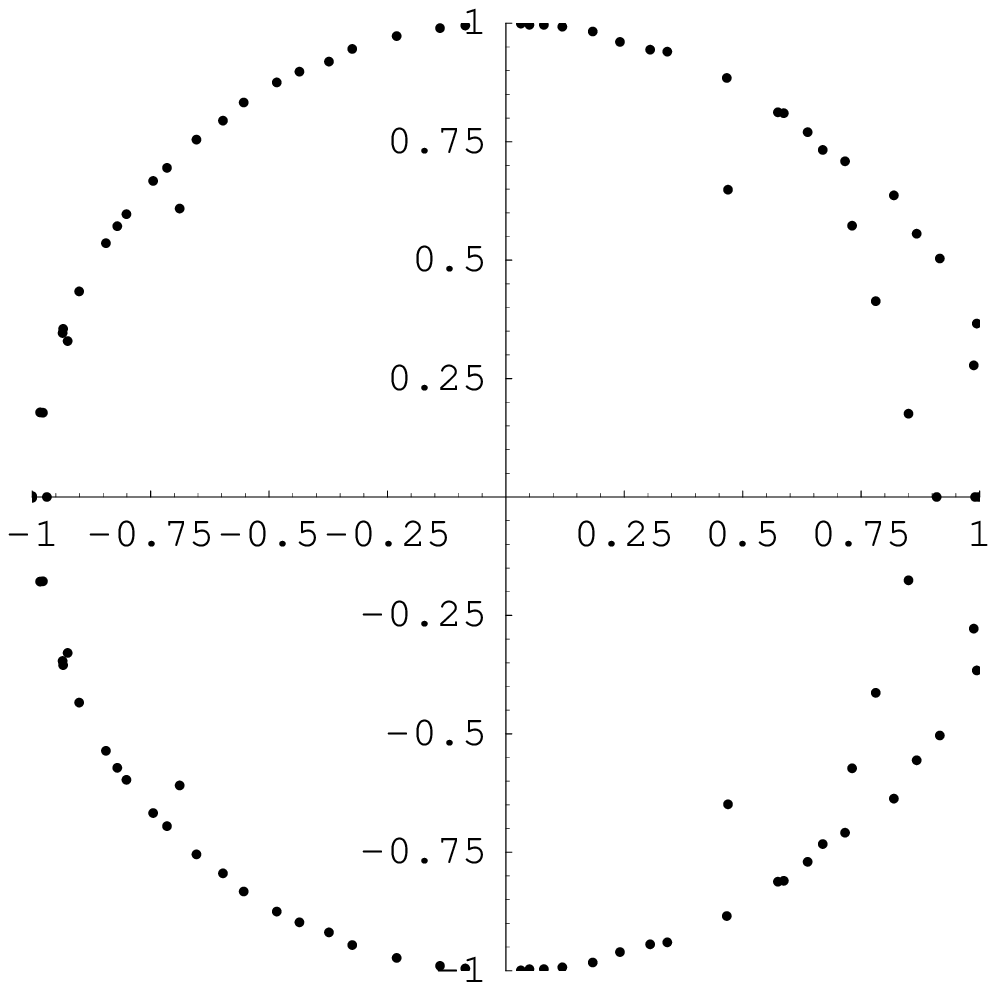}\hfill\includegraphics[scale=.50]{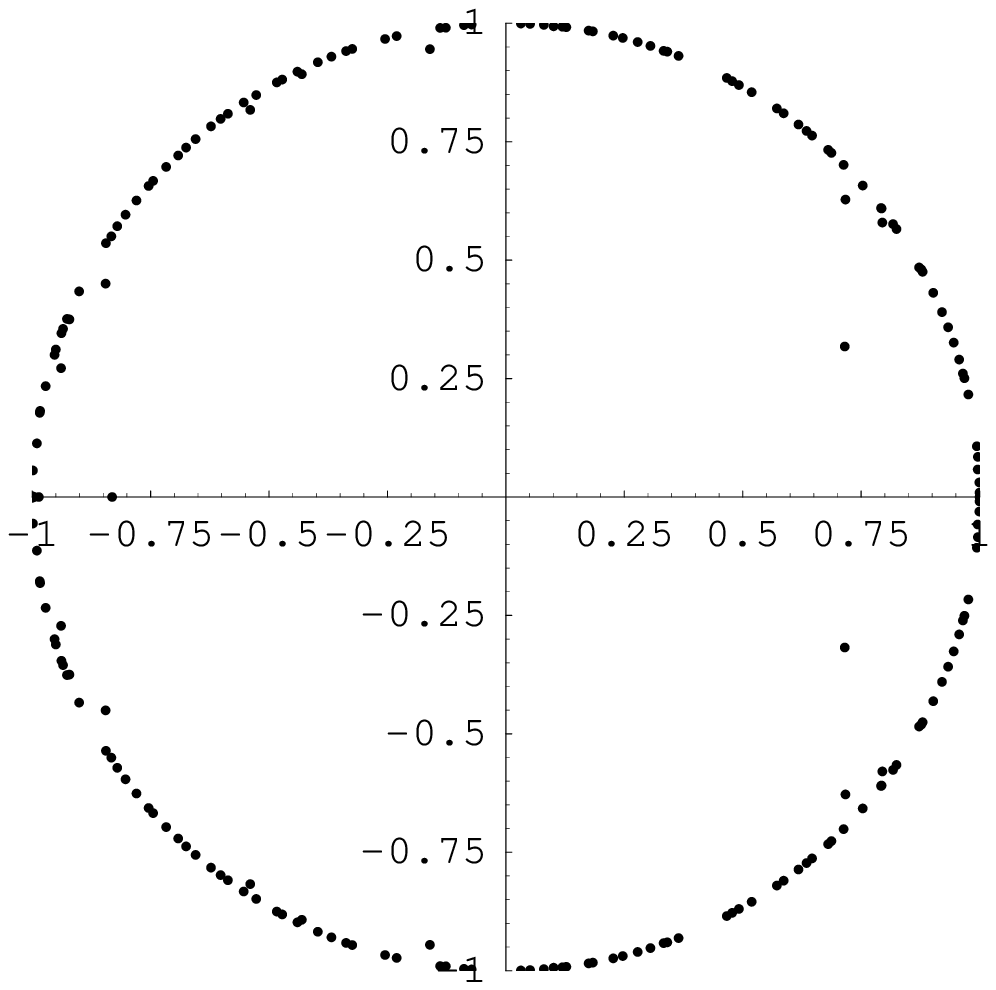}\hfill\mbox{}\\
\caption{}
\end{figure}
\end{center}

Figure~1 beautifully illustrates Theorem~\ref{T1.1} and the Nevai-Totik theorem. All the 
zeros but one clearly approach the circle $\abs{z}=\f12$. There is one zero that approaches 
approximately $0.84 + 0i$. It is the single Nevai-Totik zero in this case. That there 
are no limiting zeros inside $\abs{z}=\f12$ is no accident; see Corollary~1 of \cite{BLS}  
and Theorem~\ref{T7.2} below. And the zeros certainly seem uniformly distributed --- 
indeed, when I first ran the program that generated Figure~1, I was impressed by how 
uniform the distribution seemed to be, even for $N=10$. 

The conclusion of the Mhaskar-Saff theorem is not true for the example in Figure~2 (nor, 
of course, the hypotheses since \eqref{1.5} fails), although it would be if uniform 
density on $[-\f12,\f12]$ were replaced by any rotation invariant density, $d\gamma$,  
with $\int \abs{\log (z)}\, d\gamma(z) <\infty$ (see \cite[Theorem~10.5.19]{OPUC2}). 
But, by Theorems~10.5.19 and 10.5.21 of \cite{OPUC2}, $d\nu_n$ has a limit $d\nu$ in 
the case of Figure~2, and this limit can be seen to be of the form $f(\theta) \f{d\theta}{2\pi}$ 
with $f\in C^\infty$ and $f\not\equiv 1$ but not too far from $1$ (e.g., all odd moments 
vanish and $\int z^2 \, d\nu = (\f{1}{24})^2$). My initial thought was that the 
roughness was trying to emulate the pure point spectrum. 

I now think I was wrong in both initial reactions. 

\begin{proof}[Expectation 1: Poisson Behavior] For Figure~2, I should have made the 
connection with work of Molchanov \cite{Mol} and Minami \cite{Min} who proved in 
the case of random Schr\"odinger operators that, locally, eigenvalues in a large 
box had Poisson distribution. This  leads to a conjecture. First some notation: 

We say a collection of intervals $\Delta_1^{(n)}, \dots, \Delta_k^{(n)}$ in $\partial\bbD$ 
is canonical if $\Delta_j^{(n)} = \{e^{i\theta} \mid \theta\in [\f{2\pi a_j}{n} + \theta_j, 
\f{2\pi b_j}{n} + \theta_j]\}$ where $0\leq \theta_1 \leq \dots \leq \theta_k \leq 2\pi$,  
and if $\theta_j = \theta_{j+1}$, then $b_j < a_{j+1}$. 
\renewcommand{\qedsymbol}{}
\end{proof}

\begin{conjecture}\lb{Con2.1} Let $\{\alpha_\ell\}_{\ell=0}^\infty$ be independent identically 
distributed random variables, each uniformly distributed in $\{z\mid \abs{z}\leq \rho\}$ for 
some $\rho <1$ and let $z_1^{(n)}, \dots, z_n^{(n)}$ be the random variable describing the 
zeros of the $\Phi_n$ associated to $\alpha$. Then for some $C_1, C_2$, 
\begin{SL} 
\item[{\rm{(1)}}] 

\begin{equation} \lb{2.1} 
\bbE \bigl( \# \{z_j^{(n)} \mid \abs{z_j^{(n)}} < 1-e^{-C_1 n}\}\bigr)  \leq C_2 (\log n)^2 
\end{equation}
\item[{\rm{(2)}}] For any collection $\Delta_1^{(n)}, \dots, \Delta_k^{(n)}$ of canonical 
intervals and any $\ell_1, \dots, \ell_k$ in $\{0,1,2,\dots \}$, 
\begin{align} 
\bbP (\# \{ j \mid \arg(z_j^{(n)})\in \Delta_m\} &= \ell_m \text{ for } m=1, \dots, k) \notag \\
&\to \prod_{m=1}^k \, \biggl[ \f{(b_m -a_m)^{\ell_m}}{\ell_m!} \, e^{-(b_m -a_m)}\biggr] \lb{2.2} 
\end{align}
\end{SL} 
\end{conjecture} 

{\it Remarks.} 1. This says that, asymptotically, the distribution of $z$'s is the same 
as $n$ points picked independently, each uniformly distributed. 

\smallskip
2. See the next section for a result towards this conjecture. 

\smallskip 
3. $\bbE$ means expectation and $\bbP$ probability. 

\smallskip
4. I base the precise $e^{-C_1 n}$ and $(\log n)^2$ on the results of Stoiciu \cite{Stoi}, 
but I would regard as very interesting any result that showed, except for a small fraction 
(even if not as small as $(\log n)^2/n$), all zeros are very close (even if not as close as 
$e^{-C_1 n}$) for $\partial\bbD$. 

\smallskip 
5. There is one aspect of this conjecture that is stronger than what is proven for the 
Schr\"odinger case. The results of Molchanov \cite{Mol} and Minami \cite{Min} are the analog 
of Conjecture~\ref{Con2.1} if $\theta_1 =\theta_2 =\cdots = \theta_k$, which I would call 
a local Poisson structure. That there is independence of distant intervals is conjectured 
here but not proven in the Schr\"odinger case. That this is really an extra result can be 
seen by the fact that Figure~2 is likely showing local Poisson structure about any 
$\theta_0\neq 0,\pi$, but because the $\alpha$'s are real, the set of zeros is invariant 
under complex conjugation, so intervals about, say, $\pi/2$ and $3\pi/2$ are not independent. 

\smallskip
As far as Figure~1 is concerned, it is remarkably regular so there is an extra phenomenon 
leading to 

\begin{proof}[Expectation 2: Clock Behavior] If for $b\in (0,1)$ and $C\in \bbC$, $\alpha_n/b^n$ 
converges to $C$ sufficiently fast, then the non-Nevai-Totik zeros approach $\abs{z}=b$ 
and the angular distance between nearby zeros in $2\pi/n$. 
\renewcommand{\qedsymbol}{}
\end{proof}

\smallskip
{\it Remarks.} 1. Proving this expectation when ``sufficiently fast" means BLS convergence 
is the main new result of this paper; see Section~\ref{s4}. 

\smallskip
2. This is only claimed for local behavior. We will see that, typically, errors in the 
distance between the zeros are $O(1/n^2)$ and will add up to shift zeros that are a 
finite distance from each other relative to a strict clock. 

\smallskip
3. Clock behavior has been discussed for OPRL. Szeg\H{o} \cite{SzBk} has $C/n$ upper 
and lower bounds (different $C$'s) in many cases and Erd\"os-Turan \cite{ET40} prove 
local clock behavior under hypotheses on the measure, but their results do not cover 
all Jacobi polynomials. In Section~\ref{s6new}, we will prove a clock result for a class 
of OPRL in terms of their Jacobi matrix parameters ($\sum_{n=1}^\infty n(\abs{b_n} + 
\abs{a_n-1})<\infty$), and in Section~\ref{s6}, a simple analysis that proves local 
clock behavior for Jacobi polynomials.  I suppose this is not new, but I have not 
located a reference. 

\smallskip
4. A closer look at Figure~1 suggests that this conjecture might not be true near 
$z=b$. In fact, the angular gap there is $2(2\pi/n) +o(1/n)$, as we will see. 

\smallskip
I should emphasize that the two structures we suspect here are very different from 
what is found in the theory of random matrices. This is most easily seen by looking 
at the distribution function for distance between nearest zeros scaled to the local 
density. For the Poisson case, there is a constant density, while for clock, it is 
a point mass at a point $\theta_0 \neq 0$  depending on normalization. For the 
standard random matrix (GUE, GOE, CUE), the distribution is continuous but vanishing 
at $0$ (see \cite{MehBk}). 

Since any unitary with distinguished cyclic vector can be represented by a CMV matrix, 
CUE has a realization connected with OPUC, just not either the totally random or 
BLS case. Indeed, Killip-Nenciu \cite{KilNen} have shown that CUE is given by independent 
$\alpha_j$'s but not identically distributed. 

In Section~\ref{s3}, we describe a new result of Stoiciu \cite{Stoi} on the random case. 
In Section~\ref{s4}, we overview our various clock results: paraorthogonal OPUC in 
Section~\ref{s5}, OPRL proven in Sections~\ref{s6new} and \ref{s6}, and BLS in 
Sections~\ref{s7}--\ref{s11}. We mention some examples in Section~\ref{s11}.

\section{Stoiciu's Results on the Random Case} \lb{s3}

Recall that given $\beta\in\partial\bbD$ and $\{\Phi_n\}_{n=0}^\infty$, a set of orthogonal 
polynomials, the paraorthogonal polynomials (POPs) \cite{JNT89,Gol02} are defined by 
\begin{equation} \lb{3.1x} 
\Phi_n (z;\beta) = z\Phi_{n-1}(z) - \bar\beta \Phi_{n-1}^* (z)  
\end{equation}
They have all their zeros on $\partial\bbD$ (see, e.g., \cite[Theorem~2.2.12]{OPUC1}). 
Stoiciu \cite{Stoi} has proven the following result: 

\begin{theorem}[Stoiciu \cite{Stoi}]\lb{T3.1} Let $\{\alpha_j\}_{j=0}^\infty$ be independent 
identically distributed random variables with common distribution uniform in 
$\{z\mid\abs{z}\leq \sigma\}$ for some $\sigma <1$. Let $\{\beta_j\}_{j=0}^\infty$ be 
independent identically distributed random variables uniformly distributed on 
$\partial\bbD$. Let $z_j^{(n)}$ be the zeros of $\Phi_n (z;\beta_{n-1})$. Let 
$\Delta_1^{(n)}, \dots, \Delta_k^{(n)}$ be canonical intervals with the same 
$\theta$, that is, $\theta_1 = \theta_2 = \cdots = \theta_k$. Then \eqref{2.2} holds. 
\end{theorem} 

This differs from Conjecture~\ref{Con2.1} in two ways: The zeros are of the POPs, not 
the OPUC, and the result is only local (i.e., all $\theta$'s are equal). While the proof 
has some elements in common with the earlier work on OPRL of Molchanov \cite{Mol} and 
Minami \cite{Min}, there are many differences forced, for example, by the fact that rank  
one perturbations of selfadjoint operators differ in many ways from rank two perturbations 
of unitaries. Since the proof is involved and the earlier papers have a reputation of being 
difficult, it seems useful to summarize here the strategy of Stoiciu's proof. 

Following Minami, a key step is the proof of what are sometimes called fractional 
moment bounds and which I like to call Aizenman-Molchanov bounds after their first 
appearance in \cite{AM}. A key object in these bounds is the Green's function which 
has two natural analogs for OPUC: 
\begin{align} 
G_{nm}(z) &= \langle \delta_n, (\calC-z)^{-1} \delta_m\rangle \lb{3.1} \\ 
F_{nm}(z) &= \langle \delta_n, (\calC+z)(\calC-z)^{-1} \delta_m \rangle \lb{3.2} 
\end{align} 
These are related by 
\begin{equation} \lb{3.3} 
F_{nm}(z) =\delta_{nm} + 2z G_{nm}(z)
\end{equation}
so controlling one on $\partial\bbD$ is the same as controlling the other. 

As we will see below, $F$ and $G$ lie in the Hardy space $H^p$ for any $p<1$, so we can 
define 
\begin{equation} \lb{3.4} 
G_{nm}(e^{i\theta}) = \lim_{r\uparrow 1}\, G_{nm} (re^{i\theta}) 
\end{equation} 
for a.e.~$e^{i\theta}$. In the random case, rotation invariance will then imply that for 
any $e^{i\theta}\in\partial\bbD$, \eqref{3.4} holds for a.e.~$\alpha$. In treatments of 
Aizenman-Molchanov bounds for Schr\"odinger operators, it is traditional to prove bounds 
on the analog of $G_{ij}(z)$ for $z=re^{i\theta}$ with $r<1$ uniform in $r$ in $(\f12, 1)$. 
Instead, the Stoiciu proof deals directly with $r=1$, requiring some results on boundary  
values of $H^p$ functions to replace a uniform estimate. 

Given $N$\!, we define $\widehat\calC^{(N)}$ to be the random CMV matrix (\cite{CMV} 
and \cite[Chapter~4]{OPUC1}) obtained by setting $\alpha_N$ to the random value 
$\beta_N\in\partial\bbD$. $\widehat\calC^{(N)}$ decouples into a direct sum of an 
$N\times N$ matrix, $\calC^{(N)}$, and an infinite matrix which is identically distributed 
to the random $\calC$ if $N$ is even and $\ti\calC$, the random alternate CMV matrix, if 
$N$ is odd. (This is a slight oversimplification. Only if $\beta_N=-1$ is the infinite 
piece of $\widehat\calC^{(N)}$ a CMV matrix since the $1\times 1$ piece in the $\calM$ 
half of the $\calL\calM$ factorization has $-\beta_N$ in place of $1$. As explained in 
\cite[Theorem~4.2.9]{OPUC1}, there is a diagonal unitary equivalence taking such a 
matrix to a CMV matrix with Verblunsky coefficients $-\beta_N^{-1} \alpha_{j+N+1}$ and 
the distribution of these is identical to the distribution of the $\alpha_{j+N+1}$. 
We will ignore this subtlety in this sketch.)

We define $F_{nm}^{(N)}(z)$ and $G_{nm}^{(N)}(z)$ for $n,m\in (\{0,1,\dots, N-1\}$ 
by replacing $\calC$ in \eqref{3.1}/\eqref{3.2} by $\calC^{(N)}$. 

$\calC-\widehat\calC^{(N)}$ is a rank two matrix with $(\calC-\widehat\calC^{(N)})_{nm}\neq 0$ 
only if $\abs{n-N}\leq 2$, $\abs{m-N}\leq 2$. Moreover, any matrix element of the difference 
is bounded in absolute value by $2$. If $n\in \{0, \dots, N-1\}$, $m\geq N$\!, then 
$(\widehat\calC^{(N)}-z)_{nm}^{-1} =0$, so the second resolvent formula implies 
\begin{equation}\lb{3.4a}
\begin{aligned}  
n&\leq N-1, \, m \geq N \Rightarrow \abs{G_{nm}(z)} \\
&\leq 2 
\biggl(\, \sum_{k=N-1, N-2, N-3} \, \abs{G_{nk}^{(N)}(z)}\biggr)  
\biggl(\, \sum_{\abs{k-N+\f12}\leq \f52} \, \abs{G_{km}(z)}\biggr)
\end{aligned}
\end{equation}
which we will call the decoupling formula. 

Similarly, we have 
\begin{equation} \lb{3.4b} 
\begin{aligned}
n,m &\leq N \Rightarrow \abs{G_{nm}(z) -G_{nm}^{(N)}(z)} \\ 
&\leq  2  \biggl(\, \sum_{k=N-1, N-2, N-3} \, \abs{G_{nk}^{(N)}(z)}\biggr) 
\biggl(\, \sum_{\abs{k-N+\f12}\leq \f52}\, \abs{G_{km}(z)}\biggr)
\end{aligned}
\end{equation} 
and 
\begin{equation} \lb{3.4c} 
\begin{aligned}
n,m &\geq N \Rightarrow \abs{G_{nm}(z) -\widehat G_{nm}^{(N)}(z)} \\ 
&\leq 2 \biggl(\, \sum_{k=N,N+1,N+2} \, \abs{\widehat G_{nk}^{(N)}(z)}\biggr) 
\biggl(\, \sum_{\abs{k-N+\f12}\leq \f52}\, \abs{G_{km}(z)}\biggr)  
\end{aligned}
\end{equation}
where we recall that if $n,m\geq N$ and $N$ is even, then $\widehat G_{nm}^{(N)}
(z) = G_{n-N,m-N}  (z,\{\alpha_{j+N+1}\}_{j=0}^\infty)$ and if $N$ is odd, then 
$\widehat G_{nm}^{(N)}(z) = G_{m-N, n-N}$ $(z, \{\alpha_{j+N+1}\}_{j=0}^\infty)$. 

Stoiciu's argument has five parts, each with substeps: 
\begin{SL} 
\item[{\bf{Part 1:}}] Some preliminaries concerning $H^p$ properties of $F_{ij}$, 
positivity of the Lyapunov exponent, and exponential decay of $F_{ij}$ for a.e.~$\alpha$. 
\item[{\bf{Part 2:}}] Proof of the Aizenman-Molchanov estimates. 
\item[{\bf{Part 3:}}] Using Aizenman-Molchanov estimates to prove that eigenvalues 
of $\calC^{(N)}$ are, except for $O(\log N)$ of them, very close to eigenvalues 
of $\log N$ independent copies of $\calC^{(N/\log N)}$. 
\item[{\bf{Part 4:}}] A proof that the probability of $\calC^{(N)}$ having two 
eigenvalues in an interval of size $2\pi x/N$ is $O(x^2)$. 
\item[{\bf{Part 5:}}] Putting everything together to get the Poisson behavior. 
\end{SL} 

\smallskip
Part~2 uses Simon's formula for $G_{ij}$ (see \cite[Section~4.4]{OPUC1}) and ideas of 
Aizenman, Schenker, Friedrich, and Hundertmark \cite{newASFH}, but the details are specific 
to OPUC and exploit the rotation invariance of the distribution in an essential way. 
Part~3 uses the strategy of Molchanov-Minami with some ideas of Aizenman \cite{Aiz},  
del Rio et al.~\cite{Sim250}, and Simon \cite{Sim-aiz}. But again, there are OPUC-specific 
details that actually make the argument simpler than for OPRL. Part~4 is a new and, 
I feel, more intuitive argument than that used by Molchanov \cite{Mol} or Minami 
\cite{Min}. It depends on  rotation invariance. Part~5, following Molchanov and Minami, 
is fairly standard probability theory. Here are some of the details. 

In the arguments below, we will act as if $\log N$ and $N/\log N$ are integers rather 
than doing what a true proof does: use integral parts and wiggle blocks of size 
$[N/\log N]$ by $1$ to get $[\log N]$ of them that add to $N$\!. 

\smallskip 
\noindent{\it\ul{Step 1.1}} ($H^p$ properties of Carath\'eodory functions). A Carath\'eodory 
function is an analytic function on $\bbD$ with $F(0)=1$ and $\Real F(z)>0$. By Kolmogorov's 
theorem (see \cite[Section~4.2]{Duren}), such an $F$ is in $H^p$, $0<p<1$ with an a priori 
bound $(0<p<1$),  
\begin{equation} \lb{3.5} 
\sup_{0<r<1}\, \int \abs{F(re^{i\theta})}^p \, \f{d\theta}{2\pi} \leq \cos  
\biggl( \f{p\pi }{2}\biggr)^{-1}
\end{equation}
For any unit vector $\eta$, $\langle \eta, \f{\calC+z}{\calC-z}\, \eta\rangle$ is a 
Carath\'eodory function so, by polarization, we have the a priori bounds 
\begin{equation} \lb{3.6} 
\sup_{0<r<1}\, \int \abs{F_{nm} (re^{i\theta})}^p \, \f{d\theta}{2\pi} 
\leq 2^{2-p} \cos \biggl( \f{p\pi}{2}\biggr) 
\end{equation}
$0<p<1$, all $m,n$. The same bound holds for $F^{(N)}$. 

\smallskip 
\noindent{\it\ul{Step 1.2}} (Pointwise estimates on expectations). $H^p$ functions have 
boundary values and the mean converges (see \cite[Theorem~2.6]{Duren}), so \eqref{3.6} 
holds for $r=1$ and the $\sup$ dropped. If one averages over the random $\alpha$ 
(or random $\alpha$ and $\beta$ for $F^{(N)}$) and uses the rotation invariance to 
see that the expectation is $\theta$-independent, we find 
\begin{equation} \lb{3.7} 
\bbE (\abs{F_{nm} (e^{i\theta})}^p) \leq 2^{2-p} \cos \biggl( \f{p\pi}{2}\biggr) 
\end{equation} 
all $n,m$, all $\theta\in [0,2\pi)$, and for all $F^{(N)}$. 

\smallskip 
\noindent{\it\ul{Step 1.3}} (Positive Lyapunov exponent). By the rotation invariance 
and the Thouless formula (see \cite[Theorems~10.5.8 and 10.5.26]{OPUC2}), the density 
of zeros is $d\theta/2\pi$, the Lyapunov exponent exists for all $z$ and is given by 
(see \cite[Theorem~12.6.2]{OPUC2})
\begin{equation} \lb{3.8} 
\gamma(z) = - \tfrac12 \, \int_{\abs{z}\leq \sigma} \log (1-\abs{z}^2) \, 
\f{d^2 z}{\pi \sigma^2} + \log (\max (1,\abs{z})) 
\end{equation} 
and, in particular, $\gamma (e^{i\theta}) >0$. 

\smallskip 
\noindent{\it\ul{Step 1.4}} (Pointwise decay of $G$). Let $z_0\in\partial\bbD$ and let $\alpha$ 
be a random sequence of $\alpha$'s for which $\f{1}{n}\log \|T_n (z;\alpha)\|\to\gamma$ and 
$\abs{F_{00}(z_0)}<\infty$. Since $\gamma (z_0) >0$, the Ruelle-Osceledec theorem 
(see \cite[Theorem~10.5.29]{OPUC2}) implies there is a $\lambda\neq 1$ for which the OPUC 
with boundary condition $\lambda$ (see \cite[Section~3.2]{OPUC1}) obeys $\abs{\varphi_n^\lambda 
(z_0)}\to 0$. It follows from Theorem~10.9.3 of \cite{OPUC2} that $\abs{G_{0\,n}(z_0)}\to 0$. 
Thus, for a.e.~$\alpha$, 
\begin{equation} \lb{3.9} 
\lim_{n\to\infty} \, \abs{G_{0\,n} (z_0)}\to 0 
\end{equation} 

By Theorem~10.9.2 of \cite{OPUC2}, $F_{0\,0}\in i\bbR$ and this implies that the solutions 
$\pi$ and $\rho$ of Section~4.4 of \cite{OPUC1} obey $\abs{\pi_k (z_0)} = \abs{\rho_k(z_0)}$ 
so $\abs{(\calC-z)_{mn}^{-1}}$ is symmetric in $m$ and $n$. Thus \eqref{3.9} implies 
\[
\lim_{n\to\infty}\, \abs{G_{n\,0} (z_0)}\to 0 
\]

\smallskip 
\noindent{\it\ul{Step 1.5}} (Decay of $\bbE (\abs{G_{0\,n}(z_0)}^p)$). The proof of 
\eqref{3.9} shows, for fixed $\alpha$, $\abs{G_{0n}(z_0)}$ decays exponentially, but 
since the estimates are not uniform in $\alpha$, one cannot use this alone to conclude 
exponential decay of the expectation. But a simple functional analytic argument shows 
that if $h_n$ are functions on a probability measure space, $\sup_n \bbE (\abs{h_n}^p)
<\infty$ and $\abs{h_n(\omega)}\to 0$ for a.e.~$\omega$, then $\bbE (\abs{h_n}^r)\to 0$ 
for any $r<p$. It follows from \eqref{3.7} and \eqref{3.9} that for any $z_0\in\partial\bbD$ 
and $0<p<1$, 
\begin{equation} \lb{3.10} 
\lim_{n\to\infty}\, \bbE (\abs{G_{0\,n}(z_0)}^p) = 0 
\end{equation}

\smallskip 
\noindent{\it\ul{Step 1.6}} (General decay of $\bbE (\abs{G}^p)$). By the Schwartz 
inequality and repeated use of \eqref{3.4a}, \eqref{3.4b}, and \eqref{3.4c}, one sees 
first for $p<\f12$ and then by H\"older's inequality that 
\begin{equation} \lb{3.11} 
\lim_{n\to\infty}\, \sup_{\abs{m-k}\geq n}\, \bbE (\abs{G_{mk}(z_0)}^p) =0 
\end{equation}
and 
\begin{equation} \lb{3.12} 
\lim_{n\to\infty} \, \sup_{\substack{ \abs{m-k}\geq n \\ 0 \leq m,k\leq N-1 \\ 
\text{all } N}} \, \bbE (\abs{G_{mk}^{(N)}(z_0)}^p) =0 
\end{equation}

\smallskip
That completes Part~1. 

\bigskip 
\noindent{\it\ul{Step 2.1}} (Conditional expectation bounds on diagonal matrix elements). Let 
\begin{equation} \lb{3.13} 
H(\alpha, \gamma) = \f{\alpha + \gamma}{1+\bar\alpha \gamma} 
\end{equation} 
Then a simple argument shows that for $0<p<1$, 
\[
\sup_{\beta,\gamma} \int_{\abs{\alpha}\leq 1}\, \biggl| 
\f{1}{1-\beta H(\alpha,\gamma)}\biggr|^p\, d^2 \alpha < \infty 
\]
because, up to a constant, the worst case is $\beta =\gamma =1$, and in that case, the 
denominator vanishes only on $\Real\alpha =0$. Applying this to Khrushchev's formula  
(see \cite[Theorem~9.2.4]{OPUC2}) provides an a priori bound on the conditional 
expectation 
\begin{equation} \lb{3.14} 
\bbE (\abs{F_{kk}(z)}^p \mid \{\alpha_j\}_{j\neq k}) \leq C 
\end{equation}
where $C$ is a universal constant depending only on $\sigma$ (the radius of the 
support of the distribution of $\alpha$) and a similar result for conditioning on 
$\{\alpha_j\}_{j\neq k -1}$, that is, averaging over $\alpha_{k-1}$. 

\smallskip
\noindent{\it\ul{Step 2.2}} (Conditional expectation bounds on near diagonal matrix elements). 
Since $\rho^{-1}\leq (1-\sigma^2)^{-1/2}\equiv Q$ for all $\alpha$ with $\abs{\alpha}\leq 
\sigma$, we have, by equation~(1.5.30) of \cite{OPUC1}, that 
\[
\biggl| \f{\varphi_k}{\varphi_m}\biggr| \leq (2Q)^{\abs{k-m}} 
\]
on $\partial\bbD$. A similar estimate for the solutions $\pi$ and $\rho$ of Section~4.4 
of \cite{OPUC1} (using $\abs{\pi_k} =\abs{\rho_k}$; see the end of Step~1.4) proves  
\[ 
\biggl| \f{\pi_k}{\pi_m}\biggr| \leq (2Q)^{\abs{k-m}} 
\]
This implies, by Theorem~4.4.1 of \cite{OPUC1}, that 
\begin{equation} \lb{3.15} 
\biggl| \f{G_{kl}(z)}{G_{mn}(z)}\biggr| \leq (2Q)^{\abs{k-m} + \abs{\ell-n}} 
\end{equation} 
something clearly special to OPUC. This together with \eqref{3.14} and \eqref{3.3} 
implies 
\begin{equation} \lb{3.16} 
\bbE (\abs{G_{mn} (z)}^p\mid \{\alpha_j\}_{j\neq k}) \leq (2Q)^{\abs{m-k} + \abs{n-k}} 
(1+2C) 
\end{equation} 

\smallskip
\noindent{\it\ul{Step 2.3}} (Double decoupling). This step uses an idea of 
Aizenman, Schenker, Friedrich, and Hundertmark \cite{newASFH}. Given $n$, we look at 
$N<n-3$ and decouple first at $N$ and then at $N+3$ to get, using \eqref{3.4a} and 
\eqref{3.4c}, that 
\begin{equation} \lb{3.17} 
\begin{split}
\abs{G_{0\,n}(z)} \leq 4 \biggl(\, &\sum_{k=N-1,N-2,N-3} \, 
\abs{G_{0\,k}^{(N)}(z)}\biggr) \\
&\quad \biggl(\, \sum_{\substack{\abs{k-N+\f12}\leq \f52 \\ \abs{\ell-N+\f72} \leq \f52}}\, 
\abs{G_{k\ell}}\biggr)  \biggl(\, \sum_{\ell=N+3, N+4, N+5}\, 
\abs{\widehat G_{\ell n}^{(N+3)}(z)}\biggr) 
\end{split}
\end{equation}

Using \eqref{3.15} and generously overestimating the number of terms, we find 
\begin{equation} \lb{3.18} 
\abs{G_{0\,n}(z)} \leq 4 \cdot 3 \cdot 6 \cdot 6 \cdot 3 (2Q)^{10} 
\abs{G_{0\, N-1}^{(N)}(z)} \, \abs{G_{N+1\, N+1}(z)} \, \abs{\widehat G_{N+3\,\, n}^{(N+3)}(z)}  
\end{equation} 
Raise this to the $p$-th power and average over $\alpha_{N+1}$ with $\{\alpha_k\}_{k\neq N+1}$ 
fixed. Since $G_{0\, N-1}^{(N)}$ and $\widehat G_{N+3\, N}^{(N+3)}$ are independent of 
$\alpha_{N+1}$, they come out of the conditional expectation which can be bounded by 
\eqref{3.14}. 

After that replacement has been made, the other two factors are independent. Thus, if 
we integrate over the remaining $\alpha$'s and use the structure of $\widehat G$, we get 
\begin{equation} \lb{3.19} 
\bbE (\abs{G_{0\,n}(z)}^p) \leq C_p \bbE(\abs{G_{0\, N-1}^{(N)}(z)}^p) 
\bbE(\abs{G_{0\, n-N-3}^{(N)}}^p) 
\end{equation} 
where $C_p$ is $p$-dependent but $N$-independent. 

\smallskip
\noindent{\it\ul{Step 2.4}} (Aizenman-Molchanov bounds). By \eqref{3.12}, for $p$ 
fixed, we can pick $N$ so large that in \eqref{3.19}, we have $C_p \bbE 
(\abs{G_{0\, N-1}^{(N)}(z)}^p) <1$. If we iterate, we then get exponential decay, 
that is, we get the Aizenman-Molchanov bound; for any $p\in (0,1)$, there is 
$D_p$ and $\kappa_p$ so that 
\begin{equation} \lb{3.20} 
\bbE (\abs{G_{nm}(z)}^p) \leq D_p \, e^{-\kappa_p \abs{n-m}} 
\end{equation}
and $n,m\in [0,N-1]$, 
\begin{equation} \lb{3.21} 
\bbE (\abs{G_{nm}^{(N)}(z)}^p) \leq D_p \, e^{-\kappa_p \abs{n-m}} 
\end{equation}
One gets \eqref{3.21} from \eqref{3.20} by repeating Step~1.6. 

\smallskip 
That completes Part~2. 

\bigskip 
\noindent{\it\ul{Step 3.1}} (Dynamic localization). In the Schr\"odinger case, 
Aizenman \cite{Aiz} shows \eqref{3.21} bounds imply bounds on $\sup_t 
\abs{(e^{-ith})_{nm}}$. The analog of this has been proven by Simon \cite{Sim-aiz}. 
Thus, \eqref{3.21} implies 
\begin{equation} \lb{3.23} 
\bbE \bigl(\, \sup_\ell \abs{(\calC^\ell)_{nm}}^p\bigr) \leq 
D_p \, e^{-\kappa_p \abs{n-m}}
\end{equation} 
and similarly with $\calC$ replaced by $\calC^{(N)}$.  

\smallskip 
\noindent{\it\ul{Step 3.2}} (Pointwise a.e.~bounds). For a.e.~$\alpha$, there is 
$D(\alpha)$ so 
\begin{equation} \lb{3.34} 
[(\calC^{(N)})^\ell]_{nm} \leq D(\alpha) (N+1)^8 \, e^{-\kappa\abs{n-m}} 
\end{equation} 
for, by \eqref{3.23} and its $N$ analog with $p=\f12$, 
\begin{equation} \lb{3.35} 
\bbE \biggl(\, \sum_{\substack{ n,m,N \\ n,m\leq N}}\, \sup_\ell 
\abs{(C^\ell)_{nm}}^{1/2} (N+1)^{-4} e^{+\f12 \kappa_{1/2}\abs{n-m}}\biggr) < \infty 
\end{equation}

\smallskip
\noindent{\it\ul{Step 3.3}} (SULE for OPUC). Following del Rio, Jitomirskaya, 
Last, and Simon \cite{Sim250}, we can now prove SULE in the following form. 
For each eigenvalue $\omega_k$ of $\calC^{(N)}$, define $m_k$ to maximize the 
component of the corresponding eigenvector $u_k$ (the eigenvalues are simple), 
that is, 
\begin{equation} \lb{3.36} 
\abs{u_{k,m_k}} = \max_{\ell =1, \dots, N}\, \abs{u_{k,\ell}}  
\end{equation} 
Since 
\[ 
\f{1}{L}\, \sum_{\ell=0}^{L-1} \bar\omega_k^\ell [(\calC^{(N)})^\ell]_{nm} 
\to u_{k,n} \bar u_{k,m} 
\]
and $\max \abs{u_{k,l}} \geq N^{-1/2}$ (since $\sum_{\ell=0}^{N-1} \abs{u_{k,\ell}}^2 
=1$), \eqref{3.34} implies 
\begin{equation} \lb{3.37} 
\abs{u_{k,n}} \leq D(\alpha) (N+1)^{8.5} \, e^{-\kappa \abs{n-m_k}} 
\end{equation} 
which is what del Rio et al.~call semi-uniform localization (SULE).

\smallskip
\noindent{\it\ul{Step 3.4}} (Bound on the distribution of $u_k$). If $\abs{n-m_k} 
\geq \kappa^{-1} [(9.5)\log (N+1) + \log D(\alpha)]$, $\abs{u_{k,n}} \leq 
(N+1)^{-1}$ so $\sum_{\text{such }n} \abs{u_{k,n}}^2 \leq \f12$, and thus, 
\begin{equation} \lb{3.38} 
\abs{u_{k,m_k}}^2 \geq \f{1}{C[\log(N)+\log D(\alpha)]}  
\end{equation}
Since $\sum_{k=1}^N \abs{u_{k,\ell}}^2 =1$ for each $\ell$, \eqref{3.38} implies 
for each $m$, 
\[
\# \{k\mid m_k =m\} \leq C [\log (N)+ \log D(\alpha)] 
\]

\smallskip
\noindent{\it\ul{Step 3.5}} (Decoupling except for bad eigenvalues).  Let $(C^\sharp)^{(N)}$ 
be the matrix obtained from $\calC^{(N)}$ by decoupling in $\log (N)$ blocks of size 
$N/\log(N)$ where decoupling is done with random values of $\beta_{jN/\log(N)}$ in 
$\partial\bbD$. Call an eigenvalue of $\calC^{(N)}$ bad if its $m_k$ lies within 
$C_1 [\log(N)+\log D(\alpha)]$ and good if not. A good eigenvector is centered at an 
$m_k$ well within a single block and, by taking $C_1$ large, is of order at most 
$O(N^{-2})$ at the decoupling points. It follows, by using trial functions, that 
good eigenvalues move by at most $C_2 N^{-2}$ if $\calC^{(N)}$ is replaced by 
$(\calC^\sharp)^{(N)}$. 

\smallskip
\noindent{\it\ul{Step 3.6}} (Decoupling of probabilities). Fix the $k$ intervals 
of Theorem~\ref{T3.1}. We claim if $z_j^{(N)}$ are the eigenvalues of $\calC^{(N)}$ 
and $z_j^{\sharp(N)}$ of $\calC^{\sharp(N)}$, then 
\begin{equation} \lb{3.39} 
\begin{split}
\bbP (\# (\{j\mid &z_j^{(N)}\in \Delta_m^{(N)}\}=\ell_m,\, m=1, \dots, k) \\
&- \bbP (\# ( j\mid z_j^{\sharp(N)} \in \Delta_m^{(N)}) =\ell_m, \, m=1, \dots, k) \to 0
\end{split}
\end{equation} 
This follows if we also condition on the set where $D(\alpha)\leq D$ because 
the distribution of bad eigenvalues conditioned on $D(\alpha)\leq D$ is rotation 
invariant, and so the conditional probability is rotation invariant. Thus, with 
probability approaching $1$, no bad eigenvalues lie in the $\Delta_m^{(N)}$. 
Also, since the conditional distribution of good eigenvalues is $d\theta/2\pi$, 
they will lie within $O(N^{-2})$ of the edge with probability $N^{-1}$. Thus 
\eqref{3.39} holds with the conditioning. Since $\lim_{D\to\infty} (D(\alpha)>D)
\to 0$, \eqref{3.39} holds.  

\smallskip
That concludes Part~3 of the proof. For the fourth part, we note that the POP 
\begin{equation} \lb{3.33x}  
\Phi_n = z\Phi_{n-1} -\bar\beta_n \Phi_{n-1}^* =0 \Leftrightarrow 
\f{\beta_n z\Phi_{n-1}}{\Phi_n^*} =1
\end{equation} 

\bigskip
\noindent{\it\ul{Step 4.1}} (Definition and properties of $\eta_n$). Define 
$\eta_n (\theta;\alpha_0, \dots, \alpha_{n-1},\beta_n)$ for $\theta\in [0,2\pi]$ 
\begin{equation} \lb{3.34x} 
\f{\beta_n z\Phi_{n-1}}{\Phi_{n-1}^*} = \left. \exp (i\eta_n)\right|_{z=e^{i\theta}} 
\end{equation} 
The ambiguity in $\eta_n$ is settled by usually thinking of it as only defined mod $2\pi$, 
that is, in $\bbR/2\pi\bbZ$. $\eta_n$ is then real analytic and has a derivative 
$d\eta_n/d\theta$ lying in $\bbR$. We first claim 
\begin{equation} \lb{3.35x} 
\f{d\eta_n}{d\theta} >0 
\end{equation} 
for if $\ti\eta (\theta)$ is defined by $e^{i\ti\eta}=(z-z_0)/(1-z\bar z_0)$ for $z_0 
\in\bbR$, and $z=e^{i\theta}$, then 
\begin{equation} \lb{3.36x} 
\f{d\ti\eta}{d\theta} = \f{1-\abs{z_0}^2}{\abs{e^{i\theta}-z_0}^2} >0  
\end{equation}
from which \eqref{3.35x} follows by writing $\Phi_{n-1}$ as the product of its zeros, 
all of which lie in $\bbD$. We also have 
\begin{equation} \lb{3.37x} 
\int_0^{2\pi} \f{d\eta_n}{d\theta} \, d\theta = 2\pi n  
\end{equation}
This follows from the argument principle if we note $\beta_n z\Phi_{n-1}/\Phi_{n-1}^*$ 
is analytic in $\bbD$ with $n$ zeros there. Alternatively, since the Poisson kernel maps 
$1$ to $1$, \eqref{3.36x} implies $\int \f{d\ti\eta}{2\theta}\f{d\theta}{2\pi} =1$, which 
also yields \eqref{3.37x}. 

\smallskip
\noindent{\it\ul{Step 4.2}} (Independence of $\eta_n (e^{i\theta})$ and $\f{d\eta_n}{d\theta} 
(e^{i\theta_1})$). $\beta_n$ drops out of $d\eta_n/d\theta$ at all points. On the other hand, 
$\beta_n$ is independent of $z\Phi_{n-1}/\Phi_{n-1}^*$ and uniformly distributed. It follows 
that $\eta_n (e^{i\theta_0})$ and $\f{d\eta_n}{d\theta} (e^{i\theta_1})$ at any $\theta_0$ 
and $\theta_1$ are independent. Moreover, $\eta_n (e^{i\theta})$ is uniformly distributed. 

\smallskip
\noindent{\it\ul{Step 4.3}} (Calculation of $\bbE (d\eta_n/d\theta)$). As noted, 
$d\eta_n/d\theta$ is only dependent on $\{\alpha_j\}_{j=0}^{n-2}$ and, by rotation 
covariance of the $\alpha$'s (see \cite[Eqns.~(1.6.62)--(1.6.68)]{OPUC1}), 
\[
\f{d\eta_n}{d\theta}\, (\theta_0; e^{-(j+1)\varphi}\alpha_j) = 
\f{d\eta_n}{d\theta}\, (\theta_0-\varphi; \alpha_j) 
\]
It follows that since the $\alpha_j$'s are rotation invariant that $\bbE 
(\f{d\eta_n}{d\theta}(\theta_0))$ is independent of $\theta_0$ and so, by applying 
$\bbE$ to \eqref{3.37x},  
\begin{equation} \lb{3.38x} 
\bbE \biggl( \f{d\eta}{d\theta}\, (\theta_0)\biggr) = n  
\end{equation}

\smallskip
\noindent{\it\ul{Step 4.4}} (Bound on the conditional expectation). Let $I_n$ be an interval 
on $\partial\bbD$ of size $2\pi y/n$. Let $\lambda_0\in I_n$ and consider the conditional 
probability 
\begin{equation} \lb{3.39x} 
\bbP (I_n \text{ has 2 or more eigenvalues}\mid \lambda_0 \text{ is an eigenvalue})
\end{equation} 
(where we use ``eigenvalue" to refer to zeros of the POP since they are eigenvalues of 
a $\calC^{(N)}$). If that holds, $\eta_n (\lambda_0) =1$ and, for some $\lambda_1$ in 
$I$, $\eta(\lambda_1) =1$, so $\int_{I_n} \f{d\eta_n}{d\theta}d\theta \geq 2\pi$. Thus 
the conditional probability \eqref{3.39x} is bounded by the conditional probability 
\begin{equation} \lb{3.40} 
\bbP \biggl( \int_{I_n} \f{d\eta_n}{d\theta}\, d\theta \geq 2\pi \biggm| \eta(\lambda_0) 
=1 \biggr)  
\end{equation} 
While \eqref{3.39x} is highly dependent on the value of $\eta (\lambda_0)$, \eqref{3.40} 
is not since $d\eta_n/d\theta$ is independent of $\eta (\lambda_0)$. Thus, by Chebyshev's 
inequality,  
\begin{align} 
\eqref{3.39x} \leq \eqref{3.40} &\leq \bbP \biggl( \int_{I_n} \f{d\eta_n}{d\theta} 
\, d\theta \geq 2\pi \biggr) \notag \\
&\leq (2\pi)^{-1} \bbE \biggl( \int_{I_n} \f{d\eta_n}{d\theta}\, d\theta\biggr) \notag \\
&= \biggl( \f{2\pi y}{n}\biggr) (2\pi)^{-1}n = y \lb{3.41} 
\end{align} 
by \eqref{3.38x}. 

\smallskip
\noindent{\it\ul{Step 4.5}} (Two eigenvalue estimate). By a counting argument, 
\[
\begin{split}
\bbP &(\text{$I_n$ has exactly $m$ eigenvalues}) \\
&\qquad = \f{1}{m} \int_{\theta\in I_n} 
\bbP (\text{$I_n$ has exactly $m$ eigenvalues}\mid \eta_n(\theta)=1) \, d\kappa(\theta)
\end{split}
\]
where $d\kappa(\theta)$ is the density of eigenvalues which is $\f{n}{2\pi} d\theta$ 
by rotation invariance. Since $m\geq 2\Rightarrow \f{1}{m} \leq \f12$, we see 
\begin{align} 
\bbP (\text{$I_n$ has $2$ or more eigenvalues}) &\leq \tfrac12 \int_{I_n} 
\eqref{3.39x}\, d\kappa(\theta) \notag \\
&\leq \f{y}{2} \, \biggl( \f{n}{2\pi} \, \f{2\pi y}{n}\biggr) = \f{y^2}{2} \lb{3.42}
\end{align} 

The key is that for $y$ small, \eqref{3.42} is small compared to the probability that 
$I_n$ has at least one eigenvalue which is order $y$. This completes Part~4. 

\bigskip
\noindent{\it\ul{Step 5.1}} (Completion of the proof). It is essentially standard 
theory of Poisson processes that an estimate like \eqref{3.42} for a sum of a large 
number of independent point processes implies the limit is Poisson. The argument 
specialized to this case goes as follows. Use Step~3.6 to consider $\log N$ independent 
of POPs of degree $N/\log N$. The union of the $\Delta_m^{(N)}$ lies in a single 
interval, $\ti\Delta^{(N)}$, of size $C/N$ (here is where the $\theta_0 =\cdots = \theta_k$ 
condition is used) which is $y_N \f{2\pi}{(N/\log N)}$ with $y_N = C/\log N$. 
Thus the probability of any single POP having two zeros in $\ti\Delta^{(N)}$ is  
$O((\log N)^{-2})$. The probability of any of the $\log N$ POPs having two zeros is 
$O((\log N)^{-1})\to 0$. 

The probability of any single eigenvalue in a $\Delta_m^{(N)}$ is $O(1/\log N)$, 
so each interval is described by precisely the kind of limit where the Poisson 
distribution results. Since, except for a vanishing probability, no interval has 
eigenvalues from a POP with an eigenvalue in another, and the POPs are independent, 
we get independence of intervals. This completes our sketch of Stoiciu's proof of 
his result.

\section{Overview of Clock Theorems} \lb{s4}

The rest of this paper is devoted to proving various theorems about equal spacings 
of zeros under suitable hypotheses. In this section, we will state the main results 
and discuss the main themes in the proofs. It is easiest to begin with the case of 
POPs for OPUC: 

\begin{theorem} \lb{T4.1} Let $\{\alpha_j\}_{j=0}^\infty$ be a sequence of 
Verblunsky coefficients so that 
\begin{equation} \lb{4.1} 
\sum_{j=0}^\infty \, \abs{\alpha_j} <\infty 
\end{equation} 
and let $\{\beta_j\}_{j=1}^\infty$ be a sequence of points on $\partial\bbD$. Let 
$\{\theta_j^{(n)}\}_{j=1}^n$ be defined so $0\leq \theta_1^{(n)} \leq \cdots\leq 
\theta_n^{(n)} <2\pi$ and so that $e^{i\theta_j^{(n)}}$ are the zeros of the POPs 
\begin{equation} \lb{4.2} 
\Phi_n^{(\beta)}(z) = z\Phi_{n-1}(z) -\bar\beta_n \Phi_{n-1}^*(z) 
\end{equation}
Then {\rm{(}}with $\theta_{n+1}^{(n)}\equiv\theta_1^{(n)}+2\pi${\rm{)}} 
\begin{equation} \lb{4.3} 
\sup_{j=1, \dots, n} \, n \biggl| \theta_{j+1}^{(n)} - \theta_j^{(n)} 
-\f{2\pi}{n}\biggr| \to 0 
\end{equation} 
as $n\to\infty$. 
\end{theorem} 

The intuition behind the theorem is very simple. Szeg\H{o}'s theorem and Baxter's 
theorem imply on $\partial\bbD$ that (with $\varphi_n^{(\beta)} = \Phi_n^{(\beta)} 
/\|\Phi_{n-1}\|$) 
\begin{equation} \lb{4.4} 
\varphi_n^{(\beta)} (e^{i\theta})\sim e^{in\theta} \, \ol{D(e^{i\theta})^{-1}}\, 
-\bar\beta_n D(e^{i\theta})^{-1}
\end{equation}
and the zeros of the right side of \eqref{4.4} obey \eqref{4.3}! \eqref{4.4} holds only 
on $\partial\bbD$ and does not extend to complex $\theta$ without much stronger 
hypotheses on $\alpha$. That works since we know by other means that $\varphi_n^{(\beta)}$ 
has all its zeros on $\partial\bbD$. But when one looks at true OPUC, we will not have 
this a priori information and will need stronger hypotheses on the $\alpha$'s. 

There is a second issue connected with the $\sim$ in \eqref{4.4}. It means uniform 
convergence of the difference to zero. If $f_n$ and $g$ are uniformly close, $f_n$  
must have zeros close to the zeros of $g$, and we will have enough control on the right 
side of \eqref{4.4} to be sure that $\varphi_n^{(\beta)}$ has zeros near those of the 
right side of \eqref{4.4}. So uniform convergence  will be existence of zeros. 

A function like $f_n(x) =\sin(x) - \f{2}{n}\sin (nx)$, which has three zeros near $x=0$,  
shows uniqueness is a more difficult problem. 

There are essentially two ways to get uniqueness. One involves control over derivatives 
and/or complex analyticity which will allow uniqueness via an appeal to an intermediate  
value theorem or a use of Rouch\'e's theorem. These will each require extra hypotheses 
on the Verblunsky coefficients or Jacobi parameters. In the case of genuine OPUC where 
we already have to make strong hypotheses for existence, we will use a Rouch\'e argument. 

There is a second way to get uniqueness, namely, by counting zeros. Existence will 
imply an odd number of zeros near certain points. If we have $n$ such points and $n$ 
zeros, we will get uniqueness. This will be how we will prove Theorem~\ref{T4.1}. 
Counting will be much more subtle for OPRL because the close zeros will lie in 
$[-2,2]$ (if $a_n\to 1$ and $b_n\to 0$) and there can be zeros outside. For counting 
to work, we will need only finitely many mass points outside $[-2,2]$. This will be 
obtained via a Bargmann bound, which explains why our hypothesis in the next theorem 
is what it is: 

\begin{theorem}\lb{T4.2} Let $\{a_n\}_{n=1}^\infty$ and $\{b_n\}_{n=1}^\infty$ be 
Jacobi parameters that obey 
\begin{equation} \lb{4.4x} 
\sum_{n=1}^\infty n (\abs{a_n-1} + \abs{b_n}) <\infty 
\end{equation} 
Let $\{P_n\}_{n=0}^\infty$ be the monic orthogonal polynomials and let $\{E_j\}_{j=1}^\infty$ 
$(J<\infty)$ be the mass points of the associated measure which lie outside $[-2,2]$. 
Then for $n$ sufficiently large, $P_n(x)$ has precisely $J$ zeros outside $[-2,2]$ and 
the other $n-J$ in $[-2,2]$. Define $0<\theta_1^{(n)} < \cdots < \theta_{n-J}^{(n)} < 
\pi$ so that the zeros of $P_n(x)$ on $[-2,2]$ are exactly at $\{2\cos (\theta_\ell^{(n)}) 
\}_{\ell=1}^{n-J}$.  Then 
\begin{equation} \lb{4.5} 
\sup_{j=1, \dots, n-J-1}\, n \biggl| \theta_{j+1}^{(n)} - \theta_j^{(n)} - 
\f{2\pi}{2n} \biggr| \to 0 
\end{equation} 
\end{theorem} 

{\it Remarks.} 1. The Jacobi parameters are defined by the recursion relation 
\begin{equation} \lb{4.6} 
xP_n(x) = P_{n-1}(x) + b_{n+1} P_n(x) + a_n^2 P_{n-1}(x) 
\end{equation}
(with $P_0(x)=1$ and $P_{-1}(x)=0$). 

\smallskip
2. It is known for all Jacobi polynomials that the Jacobi parameters have 
$\abs{b_n}+ \abs{a_n-1} = O(n^{-2})$. So \eqref{4.4x} fails. In Section~\ref{s6}, we 
will provide a different argument that proves clock behavior for Jacobi polynomials.  

\smallskip 
3. We will also say something about $\abs{\theta_1}$ and $\abs{\pi-\theta_{n-J}}$, 
but the result is a little involved so we put the details in Section~\ref{s6new}. 

\smallskip 
Finally, we quote the result for OPUC obeying the BLS condition: 

\begin{theorem} \lb{T4.3} Suppose a set of Verblunsky coefficients obeys \eqref{1.9} for 
$C\in\bbC$, $C\neq 0$, $b\in (0,1)$, and $\Delta\in (0,1)$. Then the number, $J$, of 
Nevai-Totik points is finite, and for $n$ large, $\Phi_n(z)$ has $J$ zeros near these 
points. The other $n-J$ zeros can be written $\{z_j^{(n)}\}_{j=1}^{n-J}$ where $z_j^{(n)} 
\equiv \abs{z_j^{(n)}} e^{i\theta_j^{(n)}}$ with {\rm{(}}for $n$ large{\rm{)}} 
$\theta_0^{(n)} \equiv 0 <\theta_1^{(n)} < \theta_2^{(n)} < \cdots < \theta_{n-J}^{(n)} 
< 2\pi = \theta_{n-J+1}^{(n)}$. We have that 
\begin{gather} \lb{4.7} 
\sup_j \, \abs{\,\abs{z_j^{(n)}} -b\,} = O\biggl( \f{\log (n)}{n}\biggr) \\
 \sup_{j=0, \dots, n-J} \, n \biggl| \theta_{j+1}^{(n)} -\theta_j^{(n)} - 
\f{2\pi}{n} \biggr| \to 0  \lb{4.8} 
\end{gather} 
and 
\begin{equation} \lb{4.9a} 
\f{\abs{z_{j+1}^{(n)}}}{\abs{z_j^{(n)}}} = 1 + O\biggl( \f{1}{n\log n}\biggr) 
\end{equation}
\end{theorem} 

{\it Remarks.} 1. We will see that ``usually," the right side of \eqref{4.7} can be 
replaced by $O(1/n)$ and the right side of \eqref{4.9a} by $1+ O(1/n^2)$. 

\smallskip
2. Since $\theta_0 =0$ and $\theta_{n-J+1} =2\pi$ are not zeros, the angular gap between 
$z_1^{(n)}$ and $z_{n-J}^{(n)}$ is $2(2\pi/n)$.  

\smallskip 
3. \eqref{4.8} and \eqref{4.9a} imply that $\abs{z_{j+1}^{(n)} -z_j^{(n)}} \to 
\f{2\pi}{n} b$. 

\smallskip
The key to the proofs of Theorems~\ref{T4.2} and \ref{T4.3} will be careful 
asymptotics for $P_n$ and $\Phi_n$. For $P_n$, we will use well-known Jost-Szeg\H{o} 
asymptotics. For $\Phi_n$, our analysis seems to be new. 

We will also prove two refined results on the Nevai-Totik zeros, one of which has 
a clock! 

\begin{theorem}\lb{T4.4} Suppose that \eqref{1.4} holds with $b<1$. Let $z_0$ obey 
$\abs{z_0} >b$ and $D(1/\bar z_0)^{-1} =0$ {\rm{(}}i.e., $z_0$ is a Nevai-Totik 
zero{\rm{)}}. Let $z_n$ be zeros of $\varphi_n(z)$ near $z_0$ for $n$ large. 
Then for some $\veps >0$ and $n$ large, 

\begin{equation} \lb{4.11} 
\abs{z_n - z_0} \leq e^{-\veps n} 
\end{equation} 
\end{theorem} 

\smallskip 
{\it Remark.} In general, if $z_0$ is a zero of order $k$ of $\ol{D(1/\bar z)}^{-1}$, 
then there are $k$ choices of $z_n$ and all obey \eqref{4.11}. 

\begin{theorem} \lb{T4.5} Suppose that \eqref{1.9} holds for $C\in\bbC$, $C\neq 0$, 
$b\in (0,1)$, and $\Delta\in (0,1)$. There exists $\Delta_2,\Delta_3\in (0,1)$ so 
that if $b<z_0 < b\Delta_2$ is a zero of order $k$ of $\ol{D(1/\bar z)}^{-1}$, then for large 
$n$, the $k$ zeros of $\varphi_n(z)$ near $z_0$ have a clock pattern: 
\begin{equation} \lb{4.12}
z_n^{(j)} = z_0 + C_1 \biggl( \f{b}{\abs{z_0}}\biggr)^{n/k} 
\exp \biggl( -2\pi i \, \f{n}{k}\, \arg(z_0) \biggr) e^{2\pi i j/k}
+ O\biggl(\biggl( \f{b\Delta_3}{\abs{z_0}}\biggr)^{n/k}\biggr) 
\end{equation} 
for $j=0,1, \dots, k-1$ {\rm{(}}so the $k$ zeros form a $k$-fold clock{\rm{)}}. 
\end{theorem}

This completes the description of clock theorems we will prove in this article, but 
I want to mention three other situations where the pictures in \cite{OPUC1} suggest 
there are clock theorems plus a fourth situation: 

\begin{SL} 
\item[(A)] {\bf Periodic Verblunsky Coefficients.} As Figures~8.8 and 8.9 of \cite{OPUC1} 
suggest, if the Verblunsky coefficients are periodic (or converge sufficiently rapidly to 
the periodic case), the zeros are locally equally spaced, but are spaced inversely 
proportional to a local density of states. We will prove this in a future paper 
\cite{Saff3}. For earlier related results, see \cite{LPppt,Peh03,Peh04}. 

\item[(B)] {\bf Barrios-L\'opez-Saff \cite{BLS}} consider $\alpha_n$'s which are decaying 
as $b^n$ with a periodic modulation. A strong version of their consideration is that 
there is a period $p$ sequence, $c_1,c_2,\dots, c_p$, all nonzero with 
\begin{equation} \lb{4.10a} 
\alpha_n = b^n c_n +O((b\Delta)^n) 
\end{equation} 
with $0<\Delta <1$. In \cite{Saff2}, we will prove clock behavior for such $\alpha_n$'s. 
In general, there are $p$ missing points in the clock at $b\omega^j$ with $\omega = 
\exp (2\pi i/p)$ a $p$-th root of unity. Indeed, we will treat the more general 
\begin{equation} \lb{4.10b} 
\alpha_n =\sum_{\ell=1}^m c_\ell e^{in\theta_\ell} b^n + O((b\Delta)^n) 
\end{equation} 
for any $\theta_1, \dots, \theta_m\in [0,2\pi)$. 

\item[(C)] {\bf Power Regular Baxter Weights.} Figure~8.5 of \cite{OPUC1} (which shows 
zeros for $\alpha_n = (n+2)^{-2}$) suggests that if $\beta >1$ and $n^\beta \alpha_n 
\to C$ sufficiently fast, then one has a strictly clock result. By \cite{BLS}, all 
zeros approach $\partial\bbD$, and we believe that if their phases are $0\leq \theta_1 
\leq \theta_2 \leq \cdots \leq \theta_n < 2\pi$ and $\theta_{n+1} =\theta_1 + 2\pi$, 
then $\sup_j \abs{n[\theta_{j+1} -\theta_j]- 2\pi} \to C$. 

\item[(D)] {\bf Slow Power Decay.} Figures~8.6 and 8.7 of \cite{OPUC1} (which show  
$\alpha_n=(n+2)^{-1/2}$ and $(n+2)^{-1/8}$) were shown by Ed Saff at a conference 
as a warning that pictures can be misleading because they suggest there is a gap in 
the spectrum while we know that the Mhaskar-Saff theorem applies! In fact, I take 
their prediction of the gap seriously and suggest if $(n+2)^\beta \alpha_n\to C$ fast 
enough for $\beta <1$, then we have clock behavior away from $\theta =0$, that is, 
if $\theta_0$ is fixed and $\theta_j, \theta_{j+1}$ are the nearest zeros to 
$\theta_0$, then $n(\theta_{j+1} -\theta_j)\to 2\pi$, but that there is a single 
zero near $\theta =0$ with the next nearby zero $\theta'$ obeying $n\abs{\theta'}
\to\infty$. 
\end{SL} 

(C) and (D) present interesting open problems.

\smallskip
\section{Clock Theorems for POPs in Baxter's Class} \lb{s5}

In this section, we will prove Theorem~\ref{T4.1}. 

\begin{lemma} \lb{L5.1} If \eqref{4.1} holds, then 
\begin{equation} \lb{5.1} 
\sup_{e^{i\theta}\in\bbD}\, \biggl| \f{e^{i\theta} \varphi_{n-1} (e^{i\theta})} 
{\varphi_{n-1}^* (e^{i\theta})} - \f{e^{in\theta}\, \ol{D(e^{i\theta})}^{-1}} 
{D(e^{i\theta})^{-1}} \biggr| \to 0
\end{equation}
as $n\to\infty$. 
\end{lemma} 

{\it Remarks.} 1. Recall $\varphi_n^{(\beta)} = \Phi_n^{(\beta)}/\|\Phi_{n-1}\|$, that is, 
\begin{equation} \lb{5.1a} 
\varphi_n^{(\beta)} = z\varphi_{n-1} - \bar\beta\varphi_{n-1}^* 
\end{equation} 

\smallskip
2. Implicit here is the fact that $D(z)$ defined initially on $\bbD$ has a continuous 
extension to $\bar\bbD$. 

\begin{proof} Baxter's theorem (see \cite[Theorem~5.2.2]{OPUC1}) says that $D(z)$ lies 
in the Wiener algebra and, in particular, has a (unique) continuous extension to $\bar\bbD$,  
and that $\varphi_{n-1}^*(z)\to D(z)^{-1}$ uniformly on $\bar\bbD$ and, in particular, 
uniformly on $\partial\bbD$. \eqref{5.1a} plus $\varphi_{n-1}(e^{i\theta}) = 
e^{i(n-1)\theta} \,\ol{\varphi_{n-1}^*(e^{i\theta})}$ completes the proof. 
\end{proof} 

Since $D$ is nonvanishing on $\bar\bbD$ (see \cite[Theorem~5.2.2]{OPUC1}), the argument 
principle implies $D(e^{i\theta})= \abs{D(e^{i\theta})} e^{i\psi(\theta)}$ with $\psi$ 
continuous and $\psi (2\pi) = \psi(0)$. We will suppose $\psi(0)\in (-\pi,\pi]$. 

\begin{lemma} \lb{L5.2} For each $n$ and each $\eta\in [2\psi(0), 2\psi(0) + 2\pi)$, 
there are solutions, $\ti\theta_{j,\ti\eta}^{(n)}$, of 
\begin{equation} \lb{5.2} 
n\ti\theta + 2\psi (\ti\theta) = 2\pi j + \ti\eta 
\end{equation}
for $j=0,1,\dots, n-1$. We have that 
\begin{equation} \lb{5.2a} 
\sup_{\ti\eta, j=0,1,2,\dots, n-1} \, n \biggl|\, \ti\theta_{j+1,\ti\eta}^{(n)} - 
\ti\theta_{j,\ti\eta}^{(n)} - \f{2\pi}{n}\biggr| \to 0
\end{equation}
where $\ti\theta_{n,\ti\eta}^{(n)} \equiv 2\pi + \ti\theta_{0,\ti\eta}^{(n)}$. Moreover, 
for any $\veps$, there is an $N$ so for $n>N$\!, 
\begin{equation} \lb{5.3} 
\sup_{\abs{\eta -\eta'}< \f{\veps}{2}}\, n\abs{\ti\theta_{j, \ti\eta}^{(n)} - 
\ti\theta_{j,\ti\eta'}^{(n)}} \leq \veps
\end{equation} 
\end{lemma} 

\begin{proof} As $\ti\theta$ runs from $0$ to $2\pi$, the LHS of \eqref{5.2} 
runs from $2\psi(0)$ to $2\pi n + 2\psi(0)$. By continuity, \eqref{5.2} has solutions. 
If there are multiple ones, pick the one with $\ti\theta_{j,\ti\eta}^{(n)}$ as small 
as possible. 

Since $\psi$ is bounded, there is $C$ so 
\begin{equation} \lb{5.4} 
\biggl| \ti\theta_{j,\ti\eta}^{(n)} - \f{2\pi j}{n}\biggr| \leq \f{C}{n} 
\end{equation}
so subtracting \eqref{5.2} for $j+1$ from \eqref{5.2} for $j$, 
\begin{equation} \lb{5.5} 
\biggl|\, \ti\theta_{j+1,\ti\eta}^{(n)} - \ti\theta_{j,\ti\eta}^{(n)} -  
\f{2\pi}{n}\biggr| \leq \f{2}{n} \, \max_{\abs{\theta-\theta'}\leq \f{C+1}{n}} \, 
\abs{\psi(\theta) - \psi(\theta')} 
\end{equation} 
Since $\psi$ is continuous on $[0,2\pi]$, it is uniformly continuous, and thus the 
max in \eqref{5.5} goes to zero, which implies \eqref{5.2a}. 

To prove \eqref{5.3}, we note that subtracting \eqref{5.2} for $\ti\eta$ from 
\eqref{5.2} for $\ti\eta'$, 
\begin{equation} \lb{5.6} 
n \abs{\ti\theta_{j,\ti\eta}^{(n)} -\ti\theta_{j,\ti\eta'}^{(n)}} \leq 
\abs{\ti\eta - \ti\eta'} + 2 \abs{\psi(\ti\theta_{j,\ti\eta}^{(n)}) - 
\psi(\ti\theta_{j,\ti\eta'}^{(n)})}
\end{equation} 
\eqref{5.6} first implies 
\[
\abs{\ti\theta_{j,\ti\eta}^{(n)} - \ti\theta_{j,\ti\eta'}^{(n)}} \leq \f{C}{n} 
\]
and then implies \eqref{5.3} picking $N$ so 
\[
\sup_{\abs{\theta-\theta'}\leq \f{C}{n}}\, \abs{\psi(\theta) - \psi(\theta')} 
< \f{\veps}{4} 
\qedhere
\]
\end{proof} 

{\it Remark.} The proof shows that \eqref{5.3} continues to hold for any solutions 
of \eqref{5.2}. 

\begin{proof}[Proof of Theorem~\ref{T4.1}] The phase, $\zeta_n(\theta)$, of 
$\left. z\varphi_{n-1}/\varphi_{n-1}^*\right|_{z=e^{i\theta}}$ is monotone increasing 
and runs from $\zeta_n(0)$ to $2\pi n =\zeta_n(0)$ at $\theta =2\pi$ (see Step~4.1 
in Section~\ref{s3} for the monotonicity), so for any fixed $\beta_n =e^{-i\eta_n}$ 
with $\eta_n\in [\zeta_n(0), \zeta_n(0) + 2\pi)$, there are exactly $n$ solutions, 
$\theta_j^{(n)}$, $j=0,1,\dots, n-1$ of 
\begin{equation} \lb{5.7} 
\zeta_n (\theta_j^{(n)}) = 2\pi j + \eta_n 
\end{equation}
By \eqref{5.1}, $\theta_j^{(n)}$ solves \eqref{5.2} with 
\[
\eta' - \eta_n \to 0 
\]
as $n\to\infty$. By the remark after Lemma~\ref{L5.2}, that shows that \eqref{4.3} 
holds. 
\end{proof} 

I believe this result has an extension to a borderline where $\alpha_n = C_0 
n^{-1} + \text{ error}$, where the error goes to zero sufficiently fast 
$(\ell^1$ error may suffice; for applications, $\abs{\text{error}}\leq Cn^{-2}$ 
is all that is needed). The extension is on zeros away from $\theta =0$. 
I believe in this case that $D$ has an extension to $\partial\bbD\backslash\{1\}$ 
(see \cite[Section~12.1]{OPUC2}) and one has convergence there. Replacing 
uniformity should be some control of derivatives of $\varphi_n^*$ (an $O(n)$ bound). 
By the Szeg\H{o} mapping (see \cite[Section~13.1]{OPUC2}), this would provide another 
approach to Jacobi polynomials.

\section{Clock Theorems for OPRL With Bargmann Bounds} \lb{s6new} 

Our goal here is to prove Theorem~\ref{T4.2} which shows clock behavior 
for OPUC when \eqref{4.4x} holds. It is illuminating to consider two simple 
examples first: 

\begin{example}\lb{En6.1} Take $a_n\equiv 1$, $b_n \equiv 0$. Then 
\begin{equation} \lb{n6.1a} 
P_n (2\cos\theta) = c_n \, \f{\sin((n+1)\theta)}{\sin(\theta)}  
\end{equation} 
the Chebyshev polynomials of the second kind. The zeros are precisely at 
\[ 
\theta_j = \f{\pi j}{n +1} \qquad j =1,2,\dots, n 
\]
Note 
\begin{equation} \lb{n6.1b} 
\theta_{j+1} - \theta_j = \f{\pi}{n+1} = \f{\pi}{n} + O\biggl( \f{1}{n^2}\biggr) 
\end{equation} 
for $j=1, \dots, n-1$. In addition, 
\begin{equation} \lb{n6.1c} 
\theta_1 = \f{\pi}{n} + O\biggl( \f{1}{n^2}\biggr) \qquad \pi -\theta _n = 
\f{\pi}{n} + O\biggl( \f{1}{n^2}\biggr)
\end{equation} 
\qed 
\end{example} 

\begin{example}\lb{En6.2} Take $a_1 =\sqrt2$, $a_n=1$ ($n\geq 2$), $b_n=0$. Then 
\begin{equation} \lb{n6.1d} 
P_n (2\cos\theta) = c_n \cos (n\theta)  
\end{equation} 
the Chebyshev polynomials of the first kind. The zeros are precisely at 
\begin{equation} \lb{n6.1e} 
\theta_j = \f{\pi (j-\f12)}{n} \qquad j=1, \dots, n 
\end{equation} 
\eqref{n6.1b} holds in this case also but instead of \eqref{n6.1c}, we have 
\begin{equation} \lb{n6.1f} 
\theta_1 = \f{\pi}{2n} + O\biggl( \f{1}{n^2}\biggr) \qquad 
\pi -\theta_n = \f{\pi}{2n} + O\biggl( \f{1}{n^2}\biggr)  
\end{equation} 
\qed 
\end{example} 

We will make heavy use of the construction of the Jost function in this case. 
For Jacobi matrices, Jost functions go back to a variety of papers of Case 
and collaborators; see, for example, \cite{Case1,GC80}. I will follow ideas 
of Killip-Simon \cite{KS} and Damanik-Simon \cite{Jost1,DaSim2}, as 
discussed in \cite[Chapter~13]{OPUC2}. 

First, we note some basic facts about the zeros of $P_n(x)$, some only true when 
\eqref{4.4x} holds.  

\begin{proposition}\lb{Pn6.1} Let $d\mu$ be a measure on $\bbR$ whose Jacobi parameters 
obey \eqref{4.4x}. Then 
\begin{SL} 
\item[{\rm{(a)}}] $\supp(d\mu) =[-2,2]\cup\{E_j^+\}_{j=1}^{N_+} \cup \{E_j^-\}_{j=1}^{N_-}$ 
with $E_1^- <\cdots <E_{N_-}^- < -2$ and $E_1^+ > E_2^+ > \cdots > E_{N_+}^+ >2$. 
\item[{\rm{(b)}}] $N_+ <\infty$ and $N_- <\infty$. 
\item[{\rm{(c)}}] For any $n$, $P_n(x)$ has at most one zero in each $(E_j^-, E_{j+1}^-)$ 
{\rm{(}}$j=1, \dots, N_-${\rm{)}}, in each $(E_{j+1}^+, E_j^+)$ {\rm{(}}$j=1, \dots, N_+${\rm{)}}, 
and in $(E_{N_-}^-, -2)$ and $(2, E_{N_+}^+)$. 
\item[{\rm{(d)}}] For some $N_0$ and $n>N_0$, $P_n(x)$ has exactly one zero in each of 
the above intervals and all other zeros lie in $(-2,2)$. 
\end{SL} 
\end{proposition} 

\begin{proof} (a) holds because if $J_0$ is the free Jacobi matrix (the one with $a_n\equiv 
1$, $b_n\equiv 0$), then $J-J_0$ is compact. 

\smallskip
(b) This follows from the Bargmann bound for Jacobi matrices as proven by Geronimo 
\cite{Ger82,Ger88} and Hundertmark-Simon \cite{HunS}. 

\smallskip 
(c) That there is at most one zero in any interval disjoint from $\supp (d\mu)$ is a 
standard fact \cite{FrB}. 

\smallskip
(d) By a simple variational argument, using the trial functions in (1.2.61) of \cite{OPUC1},  
each $E_j^\pm$ is a limit point of zeros. This and (c) imply that each interval has a 
zero for large $N$\!. By a comparison argument, the $E_j$'s cannot be zeros for such 
$n$ and also shows that $\pm 2$ are not zeros. Since all zeros lie in $[E_1^-, E_1^+]$ 
(see \cite[Subsection~1.2.5]{OPUC1}), the other zeros lie in $[-2,2]$. 
\end{proof} 

{\it Remark.} It is possible $N_+$ and/or $N_-$ are zero, in which case the above proof 
changes slightly, for example, $E_1^-$ is replaced by $-2$. 

\smallskip
Next, we use the fact (see \cite[Theorem~13.6.5]{OPUC2}) that when \eqref{4.4x} holds, 
there is a Jost function, $u(z)$, described most simply in the variable $z$ with 
$E=z+z^{-1}$ ($z\in\bbD$ maps to $\bbC\backslash [-2,2]$ and $\partial\bbD$ is a twofold 
cover on $[-2,2]$ with $e^{i\theta}\to 2\cos\theta$). 

\begin{proposition}\lb{Pn6.2} Let $d\mu$ be a measure on $\bbR$ whose Jacobi parameters 
obey \eqref{4.4x}. There exists a function $u(z)$ on $\bar\bbD$, analytic on $\bbD$,  
continuous on $\bar\bbD$, and real on $\bar\bbD\cap\bbR$, so that 
\begin{SL} 
\item[{\rm{(a)}}] Uniformly on $[0,2\pi]$, 
\begin{equation} \lb{n6.1} 
(\sin\theta)  p_{n-1} (2\cos\theta) -\Ima (\, \ol{u(e^{i\theta})}\, e^{in\theta}) \to 0
\end{equation} 
\item[{\rm{(b)}}] The only zeros $u$ has in $\bbD$ are at those points $\beta_j^\pm\in 
\bbD$ with $\beta_j^\pm + (\beta_j^\pm)^{-1}=E_j^\pm$. 
\item[{\rm{(c)}}] The only possible zeros of $u$ and $\partial\bbD$ are at $z=\pm 1$, 
and if $u(\pm 1)=0$, then $\lim_{\theta\to 0} \, \theta^{-1} u(\pm e^{i\theta})=c$ 
exists and is nonzero. 
\end{SL} 
\end{proposition} 

\begin{proof} This is part of Theorems~13.6.4 and 13.6.5 of \cite{OPUC2}. 
\end{proof} 

\begin{proof}[Proof of Theorem~\ref{T4.2}] Write 
\begin{equation} \lb{n6.2} 
u(e^{i\theta}) = \abs{u(e^{i\theta})} e^{i\eta(\theta )}
\end{equation} 
Mod $2\pi$, $\eta$ is uniquely defined on $(0,\pi)$ and $(\pi,2\pi)$ since $u$ is 
nonvanishing there. 

Suppose first that $u(\pm 1)\neq 0$. Then $\eta$ can be chosen continuously at $\pm 1$ 
and so $\eta$ can be chosen continuously on $[0,2\pi]$ with 
\[
\eta (2\pi) -\eta(0) = 2\pi (N_+ + N_-) 
\]
by the argument principle and the fact that the number of zeros of $u$ in $\bbD$ is 
$N_+ + N_-$. 

It follows that if 
\begin{equation} \lb{n6.2a} 
g_n(\theta) =n\theta -\eta(\theta)  
\end{equation} 
that 
\[
g_n (2\pi) = 2\pi (n-N_+ - N_-) 
\]
and, in particular, 
\[ 
\Ima (\, \ol{u(e^{i\theta})} \, e^{in\theta}) 
\]
has at least $2(n-N_+ -N_-)$ zeros at points where $\ti\theta_j$ 
\begin{equation} \lb{n6.3} 
g_n (\ti\theta _j) = \pi j \qquad j = 0,1, \dots, 2(n-N_+ -N_-) -1 
\end{equation} 
Since $u$ is real, 
\begin{equation} \lb{n6.4} 
\ti\theta_{2(n-N_+ - N_-)-j} = 2\pi - \ti\theta_j 
\end{equation} 
and $\ti\theta_0 =0$, $\ti\theta_{n-N_+ -N_-} =\pi$. 

By the boundedness of $\eta$, \eqref{n6.2a}, and \eqref{n6.3}, 
\begin{equation} \lb{n6.5} 
\ti\theta_{j+1} - \ti\theta_j = \f{\pi}{n} + o\biggl( \f{1}{n}\biggr) 
\end{equation} 
uniformly in $j$. 

By \eqref{n6.1}, $\sin(\theta) p_{n-1} (2\cos\theta)$, $\theta\in [0,2\pi]$ has 
at least $n-N_+ - N_- + 1$ zeros on $[0,\pi]$ at points $\theta_j$ with $\theta_j - 
\ti\theta_j = o(1/n)$. Since $0$ and $\pi$ are zeros of $\sin\theta$ and for 
large $n$, $p_{n-1} (2\cos\theta)$ only has $n-N_+ -N_- -1$ zeros on $(0,\pi)$, 
$\{\theta_j\}$ are the only zeros. \eqref{n6.5} implies \eqref{4.5} and also 
\begin{equation} \lb{n6.6} 
\theta_1 = \f{\pi}{n} + O\biggl( \f{1}{n} \biggr) \qquad 
\theta_{n-N_+ - N_- -1} = \pi - \f{\pi}{n} + O\biggl( \f{1}{n}\biggr)
\end{equation} 

Next, we turn to the case where $u$ vanishes as $+1$ and/or $-1$. Suppose that 
$u(1)=0$, $u(-1)\neq 0$. By (a) of Proposition~\ref{Pn6.2} and the reality of 
$u$ (i.e., $\ol{u(e^{i\theta})} = u(e^{-i\theta})$), we have that $u(e^{i\theta}) 
= ic\theta + o(\theta)$ where $c\neq 0$, so 
\[
\eta(0) = \pm \f{\pi}{2} \qquad \eta(2 \pi) = \mp \f{\pi}{2} 
\qquad (\text{mod } 2\pi) 
\]
By the argument principle taking into account the zero at $z=1$, 
\[
\eta(2\pi) -\eta(0) =2\pi (N_+ + N_- + \tfrac12) 
\]
As in the regular case, \eqref{n6.5} holds, but in place of \eqref{n6.1}, 
\begin{equation} \lb{n6.7} 
\theta_1 = \f{\pi}{2n} + O\biggl( \f{1}{n}\biggr) \qquad 
\theta_{n-N_+ - N_- -1} = \pi - \f{\pi}{n} + O\biggl( \f{1}{n}\biggr) 
\end{equation} 
The other cases are similar. 
\end{proof} 

To summarize, we use a definition. 

\smallskip 
\noindent{\bf Definition.} We say $J$ has a resonance at $+2$ if and only if $u(1) =0$ 
and a resonance at $-2$ if $u(-1) =0$. 

\begin{theorem}\lb{Tn6.3} Let \eqref{4.4x} hold. Let $\theta_j^{(n)}$ be the points 
where $P_n (2\cos\theta) =0$, $0<\theta_1^{(n)} < \cdots < \theta_{n-N_+ - N_-}^{(n)} 
< \pi$. Then 
\[
\theta_1^{(n)} = \f{\pi}{n} + O\biggl( \f{1}{n^2}\biggr) 
\]
if $+2$ is not a resonance and 
\[
\theta_1^{(n)} = \f{\pi}{2n} + O\biggl( \f{1}{n^2}\biggr) 
\]
if $+2$ is a resonance. Similar results hold for $\theta_{n-N_+-N_-}^{(n)}$ with regard to 
a resonance at $-2$.  
\end{theorem}

\section{Clock Theorems for Jacobi Polynomials} \lb{s6}

The Jacobi polynomials, $P_n^{(\alpha,\beta)}(x)$, are defined \cite{Rain,SzBk} to 
be orthogonal for the weight 
\begin{equation} \lb{n7.1} 
w(x) =(1-x)^\alpha (1+x)^\beta 
\end{equation} 
on $[-1,1]$ where $\alpha,\beta >-1$ (to insure integrability of the weight). The $P_n$'s 
are normalized by 
\begin{equation} \lb{n7.2} 
P_n^{(\alpha,\beta)}(1) = \f{[\prod_{j=0}^{n-1} (1+\alpha -j)]}{n!} 
\end{equation} 
They are neither monic nor orthonormal, but it is known (\cite[Eqns.~(4.3.2) and (4.3.3)]{SzBk}) 
that $P_n^{(\alpha,\beta)}(x)$ differ from the normalized polynomials by $(n+1)^{1/2} 
c_n^{(\alpha,\beta)}$ with $0<\inf_n c_n^{(\alpha,\beta)} \leq \sup_n c_n^{(\alpha,\beta)} 
<\infty$ (not uniform in $\alpha,\beta$), and $2^n P_n^{(\alpha,\beta)}(x)$ have leading 
term $d_n^{(\alpha,\beta)}$ with a similar estimate to $c_n$. 

Our goal in this section is to prove  

\begin{theorem} \lb{Tn7.1} Fix $\alpha,\beta$. Let $\theta_j^{(n)}$, $j=1,2,\dots, n$, be 
defined by 
\begin{equation} \lb{n7.3} 
0<\theta_1^{(n)} < \cdots < \theta_n^{(n)} < \pi 
\end{equation} 
and $x_j^{(n)} = \cos(\theta_j^{(n)})$ are all the zeros of $P_n^{(\alpha,\beta)}(x)$. Then 
for each $\veps >0$, 
\begin{equation} \lb{n7.4} 
\sup_{j;\, \theta_j\in [\veps, \pi-\veps]}\, n \biggl|\theta_{j+1}^{(n)} - \theta_j^{(n)} -  
\f{\pi}{n}\biggr| \to 0
\end{equation} 
as $n\to\infty$. 
\end{theorem} 

{\it Remarks.} 1. It is not hard to see that $\theta_j\in [\veps,\pi-\veps]$ can be replaced by 
$\delta n < j < (1-\delta)n$ for each $\delta >0$. 

\smallskip
2. For restricted values of $\alpha,\beta$, this result is a special case of results of 
Erd\"os-Turan \cite{ET40}. Szeg\H{o} \cite{SzBk} has bounds of the form 
\[
Cn^{-1} < \theta_{j+1}^{(n)} - \theta_j^{(n)} < Dn^{-1} 
\]
with $C,D$ $\veps$-dependent. I have not found \eqref{n7.4}, but the proof depends on such 
well-known results in such a simple way that I'm sure it must be known! 

\begin{proof} We will depend on two classical results. The first is Darboux's formula 
(\cite[Theorem 8.21.8]{SzBk}) for the large $n$ asymptotics of $P_n^{(\alpha,\beta)}$: 
\begin{equation} \lb{n7.5} 
P_n(\cos\theta) = n^{-1/2} k(\theta)^{-1/2} \cos (n\theta + \gamma(\theta)) 
+ O(n^{-3/2}) 
\end{equation} 
where 
\begin{align} 
k(\theta) &= \pi^{-\f12} \sin\biggl( \f{\theta}{2}\biggr)^{-\alpha -\f12}  
\cos \biggl( \f{\theta}{2}\biggr)^{-\beta -\f12} \lb{n7.6} \\
\gamma(\theta) &= \tfrac12\, (\alpha +\beta + 1)\theta - (\alpha 
+\tfrac12)\, \tfrac{\pi}{2} \lb{n7.7} 
\end{align} 
and where the $O(n^{-3/2})$ is uniform in $\theta\in [\veps, \pi-\veps]$ for each fixed 
$\veps > 0$ and fixed $\alpha,\beta$. \eqref{n7.5} is just pointwise Szeg\H{o}-Jost 
asymptotics on $[-1,1]$ with explicit phase ($k$ is determined by the requirement 
that $k(\theta) w(\cos\theta) d (\cos\theta)$ must be a multiple of $d\theta$). 

The second formula we need (\cite[Eqn.~(13.8.4)]{Rain}) is 
\begin{equation} \lb{n7.8} 
\f{d}{dx}\, P_n^{(\alpha,\beta)}(x) = \tfrac12\, (1+\alpha + \beta + n) 
P_{n-1}^{(\alpha +1, \beta+1)} (x)
\end{equation} 
which is a simple consequence of the Rodrigues formula. 

Fix $\theta_0$. Define $\theta^{(n)}$ by 
\begin{equation} \lb{n7.9} 
\theta^{(n)} = \f{2\pi}{n}\, \biggl[ \f{n\theta_0}{2\pi}\biggr] 
\end{equation} 
where $[y]=$ integral part of $y$. Then $\theta^{(n)} \leq \theta_0 
< \theta^{(n)} +\f{2\pi}{n}$ and $n\theta^{(n)} \in 2\pi\bbZ$. Define 
\begin{equation} \lb{n7.10} 
f_n(y) =n^{1/2} P_n^{(\alpha,\beta)} \biggl( \cos\biggl( \theta^{(n)} + 
\f{y}{n}\biggr)\biggr) 
\end{equation} 
Then \eqref{n7.5} implies that uniformly in $\theta_0\in [\veps,\pi-\veps]$ and $y\in [-Y,Y]$ 
(any $\veps >0$, $Y<\infty$), we have as $n\to\infty$, 
\begin{equation} \lb{n7.11} 
f_n(y) \to k(\theta_0)^{-1/2} \cos (y+\gamma (\theta_0)) 
\end{equation} 
Moreover, by \eqref{n7.8}, $f'_n(y)$ converges uniformly to a continuous limit. Standard 
functional analysis (essentially the fundamental theorem of calculus!) says that if $f_n 
\to f_\infty$ and $f'_n\to g_\infty$ uniformly for $C^1$ functions $f_n$, then $g_\infty 
=f'_\infty$. Thus 
\begin{equation} \lb{n7.12} 
f_n' (y) \to -k(\theta_0)^{-1/2} \sin (y+\gamma(\theta_0)) 
\end{equation}

In particular, since $k(\theta_0)$ is bounded above and below on $[-\veps,\pi-\veps]$, 
we see that for $\alpha,\beta,\veps,Y$ fixed, there are $C_1,C_2 >0$ and $N$\!, so for 
all $\theta_0\in [\veps,\pi-\veps]$, $\abs{y_j} <Y$ for $j=1,2$, and $n>N$, we have 
\begin{equation} \lb{n7.13} 
\abs{y_1 - y_2} <\pi, \quad \abs{f_n(y_j)} < C_1 \Rightarrow 
\abs{f_n (y_1) - f_n(y_2)} \geq C_2 (y_1 -y_2)  
\end{equation} 

In the usual way \eqref{n7.13} implies that for $n$ large, there is at most one solution 
of $f_n(y)=0$ within $\pi$ of another solution. Since \eqref{n7.5} implies existence, 
we can pinpoint the zeros of $P_n$ in $[\veps, \pi-\veps]$ as single points near 
$\{\f{\pi}{2n} + \f{j\pi}{n} -\gamma (\f{\pi}{2n} + \f{j\pi}{n}) + o(\f{1}{n})\}$. 
From this, \eqref{n7.4} follows. 
\end{proof}

\smallskip
\section{Asymptotics Away From the Critical Region} \lb{s7}  

This is the first of several sections which focus on proving Theorem~\ref{T4.3}. 
The key will be asymptotics of $\varphi_n(z)$ in the region near $\abs{z}=b$. 
In this section, for background, we discuss asymptotics away from $\abs{z}=b$. 
We start with $\abs{z}>b$. The first part of the following is a translation of \
results of Nevai-Totik \cite{NT89} from asymptotics of $\varphi_n^*$ to $\varphi_n$: 

\begin{theorem}\lb{T7.1} If \eqref{1.4} holds for $0<b<1$, then $D$ is analytic 
in $\{z\mid \abs{z}< b^{-1}\}$ and for $\abs{z} >b$, 
\begin{equation} \lb{7.1} 
\lim_{n\to\infty} \, z^{-n} \varphi_n(z) = \ol{D(1/\bar z)}^{-1} 
\end{equation} 
and \eqref{7.1} holds uniformly in any region $\abs{z} \geq b+\veps$ with $\veps >0$. 
Indeed, on any region $\{z\mid b+\veps <\abs{z} \leq 1\}$, 
\begin{equation} \lb{7.2} 
\abs{\varphi_n(z)- \ol{D(1/\bar z)}^{-1} z^n} \leq C_\veps \biggl( b + \f{\veps}{2}\biggr)^n 
\end{equation} 
\end{theorem} 

{\it Remark.} The point of \eqref{7.2} is that the error in \eqref{7.1} is approximately 
$O(b^n/\abs{z}^n)$, which is exponentially small if $\abs{z} >b+\veps$. It is remarkable 
that we get exponentially small errors with only \eqref{1.4}. 

\begin{proof} By step (2) in the proof of Theorem~\ref{T1.1}, 
\begin{equation} \lb{7.3} 
\abs{\varphi_{n+1}^*(z) - \varphi_n^*(z)} \leq \ti C_\veps [\max (1,\abs{z})]^n 
\biggl| b + \f{\veps}{2}\, \biggr|^n 
\end{equation} 
As noted there, this implies \eqref{1.8} which, using 
\begin{equation} \lb{7.3a} 
\varphi_n(z) =z^n \, \ol{\varphi_n^* (1/\bar z)} 
\end{equation} 
implies \eqref{7.1}. \eqref{7.3} then implies 
\begin{equation} \lb{7.4} 
\abs{\varphi_n^*(z) - D(z)^{-1}} \leq C_\veps [\max(1,\abs{z})]^n \, 
\biggl| b+\f{\veps}{2} \biggr|^n
\end{equation} 
if $\abs{z} ( b+\f{\veps}{3}) <1$, and this yields \eqref{7.2} after using \eqref{7.3a}. 
\end{proof} 

{\it Remark.} The restriction $\abs{z}\leq 1$ for \eqref{7.2} comes from $\abs{z} \leq 
1$ in \eqref{1.7}. But, by Theorem~\ref{T7.1}, \eqref{1.7} holds if $\abs{z} >b$, and 
so we can conclude \eqref{7.2} in any region $\{z\mid b+\veps <\abs{z} < b^{-1}-\veps \}$.  
By the maximum principle, we have that 
\[
\sup_{\abs{z} \leq b+\veps}\, \abs{b+\veps}^{-n} \abs{\Phi_n (z)} <\infty 
\]
which, plugged into the machine in \eqref{7.2}, implies for $\abs{z} >b^{-1}-\veps$, we 
have 
\[
\abs{\varphi_n(z) - \ol{D(1/\bar z)}^{-1} z^n} \leq C_\veps \abs{z}^n 
\biggl( b+\f{\veps}{2}\biggr)^n 
\]
which is exponentially small compared to $\abs{z}^n$. 

\smallskip

Barrios-L\'opez-Saff \cite{BLS} proved that ratio asymptotics \eqref{1.8x} implies 
that $\varphi_{n+1}(z)/\varphi_n(z)\to b$ for $\abs{z} <b$, thereby also proving 
there are no zeros of $\varphi_n$ in each disk $\{z\mid \abs{z} <b-\veps\}$ if 
$n$ is large (see also \cite[Section~9.1]{OPUC2}). Here we will get a stronger 
result from a stronger hypothesis: 

\begin{theorem}\lb{T7.2} Suppose that 
\begin{equation} \lb{7.5} 
b^{-n} \alpha_n \to C\neq 0 
\end{equation}
as $n\to\infty$. Then for any $\abs{z}<b$, 
\begin{equation} \lb{7.6} 
b^{-n} \varphi_n (z) \to \bar C(z-b)^{-1} D(z)^{-1} 
\end{equation} 
Moreover, if BLS asymptotics \eqref{1.9} holds, then in each region $\{z\mid\abs{z}  
< b-\veps\}$, 
\begin{equation} \lb{7.6a} 
\abs{b^{-n} \varphi_n(z) - \bar C(z-b)^{-1} D(z)^{-1}} \leq C_1 \ti\Delta^n 
\end{equation} 
for some $\ti\Delta <1$. 
\end{theorem} 

{\it Remark.} In \cite{Saff2}, we will prove a variant of this result that 
only needs ratio asymptotics as an assumption. 

\begin{proof} Define 
\begin{equation} \lb{7.7} 
u_n(z)= \varphi_n(z) b^{-n} \qquad A_n =-\bar\alpha_n b^{-n-1} 
\end{equation} 
Then, Szeg\H{o} recursion says 
\begin{equation} \lb{7.8} 
u_{n+1} = \biggl( \f{z}{b}\biggr) u_n + A_n \varphi_n^*(z) 
\end{equation} 
Iterating, we see that 
\begin{equation} \lb{7.9} 
u_n = \sum_{j=1}^n A_{n-j} \varphi_{n-j}^* (z) \biggl( \f{z}{b}\biggr)^{j-1} + 
\biggl( \f{z}{b}\biggr)^n u_0 
\end{equation} 

Since $A_m\to -\bar Cb^{-1}$, $\varphi_m^*\to D^{-1}$, and $\abs{z}/b <1$, \eqref{7.9} implies 
$u_n$ has a limit $u_\infty$. \eqref{7.8} then implies 
\begin{equation} \lb{7.10} 
bu_\infty = zu_\infty -\bar CD(z)^{-1} 
\end{equation} 
which implies \eqref{7.6}. 

If \eqref{1.9} holds, then 
\begin{equation} \lb{7.11} 
A_n -A_\infty = O(\Delta^n) 
\end{equation} 
Moreover, Szeg\H{o} recursion for $\Phi_n^*$ implies if $\abs{z} <1$, 
\[
\abs{\Phi_{n+1}^* - \Phi_n^*} \leq C_1 b^n \abs{\rho_n^{-1} -1} 
\]
and then since $\abs{\rho_n^{-1} -1} \leq C_1 \abs{\alpha_n}$ if $\abs{\alpha_n} 
<\f14$, $\abs{\varphi_{n+1}^* - \varphi_n^*} \leq C_2 b^n$, and so
\begin{equation} \lb{7.12} 
\abs{z}\leq 1 \Rightarrow \abs{\varphi_n^*(z) - D(z)^{-1}} \leq C_3 b^n 
\end{equation} 
\eqref{7.9}, \eqref{7.11}, and \eqref{7.12} imply \eqref{7.6a} with $\ti\Delta = 
\max (\Delta, b)$. 
\end{proof}

\section{Asymptotics in the Critical Region} \lb{s8}

The key result in controlling the zeros when the BLS condition holds is 

\begin{theorem} \lb{T8.1} Let the BLS condition \eqref{1.9} hold for a 
sequence, $\{\alpha_n\}_{n=0}^\infty$, of Verblunsky coefficients and 
some $b\in (0,1)$ and $C\in\bbC$. Then there exist $D$, $\Delta_1$, 
and $\Delta_2$ with $0<\Delta_1 < \Delta_2 <1$ so that if 
\begin{equation} \lb{8.1} 
b\Delta_2 < \abs{z} < b\Delta_2^{-1}  
\end{equation}
then 
\begin{equation} \lb{8.2} 
\abs{\varphi_n(z) -\ol{D(1/\bar z)}^{-1} z^n - \bar C(z-b)^{-1} D(z)^{-1} b^n} 
\leq D(b\Delta_1)^n  
\end{equation} 
\end{theorem} 

{\it Remarks.} 1. Implicit in \eqref{8.2} is that $\ol{D(1/\bar z)}^{-1}$ has an 
analytic continuation (except at $z=b$; see Remark~2) to the region \eqref{8.1}, 
that is, $D(z)^{-1}$ has an analytic continuation to $\{z\mid\abs{z}<b^{-1} 
\Delta_2^{-1}\}$ except for $z=b^{-1}$.  

\smallskip 
2. Since $\varphi_n(z)$ is analytic at $z=b$ and $D(b)^{-1}\neq 0$, the poles in 
$\ol{D(1/\bar z)}^{-1}$ and $D(z)/(z-b)^{-1}$ must cancel, that is, \eqref{1.3} 
must hold. 

\smallskip
3. In this way, Theorem~\ref{T8.1} includes a new proof of one direction of 
Theorem~\ref{T1.2}, that is, that the BLS condition implies that $D(z)^{-1}$ 
is meromorphic in $\{z\mid\abs{z}<b^{-1}\Delta_2^{-1}\}$ with a pole only at 
$z=b^{-1}$ with \eqref{1.3}. 

\smallskip 
4. The condition $\Delta_1 < \Delta_2$ implies that the error $O((b\Delta_1)^n)$ 
is exponentially smaller than both $z^n$ and $b^n$ in the region where \eqref{8.1} 
holds. 

\smallskip 
We will prove \eqref{8.2} by considering the second-order equation obeyed by 
$\varphi_n$ for $n$ so large that $\alpha_{n-1}\neq 0$ (see \cite[Eqn.~(1.5.47)]{OPUC1}): 
\begin{equation} \lb{8.3} 
\varphi_{n+1} = \rho_n^{-1} \biggl( z + \f{\bar\alpha_n}{\bar\alpha_{n-1}}\biggr) 
\varphi_n - \f{\bar\alpha_n}{\bar\alpha_{n-1}}\, \f{\rho_{n-1}}{\rho_n}\, z\varphi_{n-1}
\end{equation} 
(the only other applications I know of this formula are in \cite{BLS} and 
Mazel et al.~\cite{MGH}). By \eqref{1.9} which implies $\rho_n=1 + O(b^{2n})$ 
and $\bar\alpha_n/\bar\alpha_{n-1} = b+O(\Delta^n)$, we have 
\begin{align} 
\rho_n^{-1} \biggl( z + \f{\bar\alpha_n}{\bar\alpha_{n-1}}\biggr) 
&= z+b+O(b^{2n} + \Delta^n) \lb{8.4} \\
\biggl( \f{\bar\alpha_n}{\bar\alpha_{n-1}}\biggr) \f{\rho_{n-1}}{\rho_n}\, z 
&= bz+ O(b^{2n} + \Delta^n) \lb{8.5} 
\end{align} 

In \cite{Saff2}, we will analyze this critical region by an alternate method 
that, instead of analyzing \eqref{8.3} as a second-order homogeneous difference 
equations, analyzes the more usual Szeg\H{o} recursion as a first-order inhomogeneous 
equation. 

We thus study the pair of difference equations: 
\begin{align}
u_{n+1} &= \rho_n^{-1} \biggl( z + \f{\bar\alpha_n}{\bar\alpha_{n-1}}\biggr) u_n - 
\f{\bar\alpha_n}{\bar\alpha_{n-1}}\, \f{\rho_{n-1}}{\rho_n}\, zu_{n-1}  \lb{8.6} \\
\intertext{and} 
u_{n+1}^{(0)} &= (z+b) u_n^{(0)} - bz u_{n-1}^{(0)} \lb{8.7} 
\end{align} 
expanding solutions of $u_n$ in terms of solutions of $u_n^{(0)}$. 

Two solutions of \eqref{8.7} are $b^n$ and $z^n$ (since $x^2 - 
(z+b)x + bz$ is solved by $x=b,z$). These are linearly independent if $b\neq z$ but not at 
$b=z$, so it is better to define 
\begin{equation} \lb{8.8} 
x_n = z^n \qquad y_n = \f{z^n - b^n}{z-b}  
\end{equation} 
with $y_n$ interpreted as $nb^{n-1}$ at $z=b$. 

We rewrite \eqref{8.7} as 
\begin{equation} \lb{8.9} 
\binom{u_{n+1}^{(0)}}{u_n^{(0)}} = M^{(0)} \binom{u_n^{(0)}}{u_{n-1}^{(0)}}  
\end{equation}
with 
\begin{equation} \lb{8.10} 
M^{(0)} = \begin{pmatrix} 
z+b & -bz \\ 1 & 0 \end{pmatrix} 
\end{equation} 
and \eqref{8.6} as
\begin{equation} \lb{8.11} 
\binom{u_{n+1}}{u_n} = (M^{(0)} + \delta M_n) \binom{u_n}{u_{n-1}}  
\end{equation} 
where $\delta M_n$ is affine in $z$ and, by \eqref{8.4}/\eqref{8.5}, obeys 
\begin{equation} \lb{8.12} 
\|\delta M_n\| \leq C (1+\abs{z}) (b^{2n} + \Delta^n)
\end{equation}

Now we use variation of parameters, that is, we define $c_n,d_n$ by 
\begin{align} 
u_n &= c_n x^n + d_n y^n  \lb{8.13} \\
u_{n-1} &= c_n x^{n-1} + d_n y^{n-1} \lb{8.14} 
\end{align} 
or $(u_n u_{n-1})^t = Q_n (c_n d_n)^t$ where 
\begin{equation} \lb{8.15} 
Q_n = \begin{pmatrix} 
x^n & y^n \\ x^{n-1} & y^{n-1}\end{pmatrix} 
\end{equation}
Since $\det(Q_n)=-z^{n-1} b^{n-1}$, we have 
\begin{equation} \lb{8.16} 
Q_{n+1}^{-1} = -z^{-n} b^{-n} \begin{pmatrix} 
y^n & -y^{n+1} \\ -x^n & x^{n+1} \end{pmatrix} 
\end{equation} 

Since $x_n$ and $y_n$ solve \eqref{8.7}, 
\begin{equation} \lb{8.17} 
Q_{n+1}^{-1} M^{(0)} Q_n = \boldsymbol{1} 
\end{equation} 
Moreover, since 
\begin{equation} \lb{8.18} 
\abs{x_n} \leq \abs{z}^n \qquad \abs{y_n}\leq n \,\,\, \max(\abs{z},\abs{b})^n 
\end{equation} 
\eqref{8.12}, \eqref{8.15}, and \eqref{8.16} imply that 
\begin{equation} \lb{8.19} 
\delta\wti M_n \equiv Q_{n+1}^{-1} \delta M_n Q_n 
\end{equation} 
obeys
\begin{equation} \lb{8.20} 
\|\delta\wti M_n\| \leq C_n \biggl[ \f{\max(\abs{z},\abs{b})}{\min(\abs{z},\abs{b})} 
\biggr]^n (1+\abs{z}) (b^{2n} + \Delta^n) 
\end{equation} 

In particular, in the region \eqref{8.1}, 
\begin{equation} \lb{8.21} 
\|\delta\wti M_n\| \leq C \Delta_1^{2n} 
\end{equation} 
if we take $\Delta_1 = \Delta_2^2$ and $\Delta_2< 1$ is picked 
so that $\max (b^2,\Delta) < \Delta_2$. Since, in the region \eqref{8.1}, 
\begin{equation} \lb{8.22} 
\abs{z^n} \leq \abs{b}^n \Delta_2^{-n} 
\end{equation} 
we have that 
\begin{equation} \lb{8.23} 
n [\abs{z}^n + \abs{b}^n ] \|\delta\wti M_n\| \leq C (b\Delta_1)^n  
\end{equation} 

We thus have the tools to prove the main input needed for Theorem~\ref{T8.1}: 

\begin{proposition} \lb{P8.2} There exist $0<\Delta_1 < \Delta_2 <1$ and 
$N$ so that for $z$ in the region \eqref{8.1}, there are two solutions $u_n^+(z)$ 
and $u_n^-(z)$ of \eqref{8.6} for $n\geq N$ with 
\begin{SL} 
\item[{\rm{(i)}}] 
\begin{equation} \lb{8.24} 
\abs{u_n^+ - x^n} - \abs{u_n^- -y^n} \leq D_1 (b\Delta_1)^n  
\end{equation} 
\item[{\rm{(ii)}}] $u_n^\pm (z)$ are analytic in the region \eqref{8.1}. 
\item[{\rm{(iii)}}] $(u_{n+1}^\pm, u_n^\pm)^t$ are independent for $^+$ and 
$^-$ for $n\geq N$\!. 
\end{SL} 
\end{proposition} 

\begin{proof} If $u_n$ is related to $c_n,d_n$ by \eqref{8.13}/\eqref{8.14}, then 
\eqref{8.6} is equivalent to \eqref{8.11} and then to 
\begin{equation} \lb{8.25} 
\binom{c_{n+1}}{d_{n+1}} = (1+\delta\wti M_n) \binom{c_n}{d_n} 
\end{equation}
By \eqref{8.23} and the fact that $\delta\wti M_n$ is analytic in the region \eqref{8.2}, 
we see 
\begin{equation} \lb{8.26} 
L_n(z) = \prod_{j=n}^\infty (1+\delta\wti M_n(z))  
\end{equation} 
exists, is analytic in $z$ (in the region \eqref{8.1}), and invertible for all $z$ 
in the region and $n\geq N$ sufficiently large.

Define 
\begin{equation} \lb{8.27} 
\binom{c_n^+}{d_n^+} = L_n(z)^{-1} \binom{1}{0}  
\end{equation} 
and $c_n^-,d_n^-$ by the same formula with $\binom{1}{0}$ replaced by 
$\binom{0}{1}$ and then $u_n^\pm$ by \eqref{8.13}. Analyticity of $u$ is immediate 
from the analyticity of $L_n(z)$ and \eqref{8.24} follows from \eqref{8.23}. 

Independence follows from the invertibility of $L_n(z)^{-1}$. 
\end{proof} 

\begin{proof}[Proof of Theorem~\ref{T8.1}] By the independence, we can write 
\begin{equation} \lb{8.28} 
(\varphi_{N+1}(z), \varphi_N(z)) = f_1(z) (u_{N+1}^+(z), u_N^+(z)) + 
f_2(z) (u_{N+1}^-(z), u_N^-(z)) 
\end{equation} 
where $f_1,f_2$ are analytic in the region \eqref{8.1} since $\varphi$ and $u^\pm$ 
are. By the fact that $\varphi,u^\pm$ obey \eqref{8.6}, we have for all $n\geq N$ that 
\begin{equation} \lb{8.29} 
\varphi_n(z) = f_1(z) u_n^+(z) + f_2 (z) u_n^- (z) 
\end{equation}

Suppose first $\abs{z}<b$ in the region \eqref{8.1}. Then, since $b^{-n}\abs{z}^n 
\to 0$, \eqref{8.24} implies that  as $n\to\infty$, 
\begin{equation} \lb{8.30} 
b^{-n} u_n^+ \to 0 \qquad b^{-n} u_n^- \to -\f{1}{z-b} 
\end{equation} 
We conclude, by \eqref{7.6}, that in that region, 
\begin{equation} \lb{8.31} 
f_2 (z) = -\bar C D(z)^{-1} 
\end{equation} 
so, by analyticity, this holds in all of the region \eqref{8.1}. 

Next, suppose $\abs{z}>b$, so $\abs{z}^{-n} b^n\to 0$. \eqref{8.24} 
implies that as $n\to\infty$, 
\begin{equation} \lb{8.32} 
z^{-n} u_n^+ \to 1 \qquad z^{-n} u_n^- \to \f{1}{z-b} 
\end{equation} 
By \eqref{7.1}, we conclude that 
\begin{equation} \lb{8.33} 
f_1(z) + f_2(z) (z-b)^{-1} = \ol{D(1/\bar z)}^{-1}  
\end{equation}
Again, by analyticity, this holds in all of the region \eqref{8.1} except $z=b$. 
It follows that $D(z)^{-1}$ has a pole at $b^{-1}$ and otherwise is analytic 
in $\{z\mid\abs{z}<b^{-1} \Delta_2^{-1}\}$. \eqref{8.33} also determines the 
residue to be given by \eqref{1.10}. 

\eqref{8.29} becomes 
\begin{equation} \lb{8.34} 
\varphi_n(z) = [\,\ol{D(1/\bar z)}^{-1} + \bar C(z-b)^{-1} D(z)^{-1}] 
u_n^+(z) - \bar CD(z)^{-1} u_n^-(z) 
\end{equation} 
\eqref{8.2} follows from this result, boundedness of $f_1,f_2$, and \eqref{8.24}. 
\end{proof}

\section{Asymptotics of the Nevai-Totik Zeros} \lb{s10new} 

In this section, we prove Theorems~\ref{T4.4} and \ref{T4.5}. They will be simple  
consequences of Theorems~\ref{T7.1} and \ref{T8.1}. 

\begin{proof}[Proof of Theorem~\ref{T4.4}] If $\ol{D(1/\bar z)}^{-1}$ has a zero of 
order exactly $k$, at $z_0$ there is some $C_1 >0$ and $\delta >0$ so that 
\begin{equation} \lb{n10.1} 
\abs{z-z_0} <\delta \Rightarrow \abs{D(1/\bar z)} \geq C_1 \abs{z-z_0}^k
\end{equation} 
For $n$ large, Hurwitz's theorem implies $\abs{z_n -z_0}<\delta$. By \eqref{7.2} 
and \eqref{n10.1}, we have 
\[
C_1 \abs{z_n-z_0}^k \leq C_d z_n^{-n} \biggl( b+\f{d}{2}\biggr)^n 
\]
for small $d$ and $n$ large. Picking $d$ so $z_0^{-1} (b+\f{d}{2})<1$, we see 
\[
\abs{z_n-z_0}^k \leq C_2 e^{-2k\veps n} 
\]
for some $\veps >0$. This implies \eqref{4.11} for $n$ large. 
\end{proof} 

\begin{proof}[Proof of Theorem~\ref{T4.5}] Pick $\Delta_2$ to be given by 
Proposition~\ref{P8.2}. Define $Q(z)$ near $z_0$ by 
\begin{equation} \lb{n10.2} 
(z_n-z_0)^k Q(z) = \ol{D(1/\bar z)}^{-1}  
\end{equation} 
so $Q(z_0) \neq 0$ and $Q$ is analytic near $z_0$. $\varphi_n (z_n) =0$ 
and \eqref{8.2} implies 
\begin{equation} \lb{n10.3} 
(z_n -z_0)^k Q(z_n) = C (z_n -b)^{-1} D(z_n)^{-1} \, \f{b^n}{z_n^n} + 
O\biggl( \f{(b\Delta_1)^n}{z_n^n}\biggr) 
\end{equation} 
By Theorem~\ref{T4.4}, $z_n -z_0= O(e^{-\veps n})$ so \eqref{n10.3} becomes 
\begin{equation} \lb{n10.4} 
(z_n -z_0)^k = C_1 \, \f{b^n}{z_0^n} + O\biggl( \f{(b\Delta_1)^n}{z_0^n}\biggr) 
+ O\biggl( \f{b^n}{z_0^n}\, e^{-\veps n}\biggr) 
\end{equation} 
since $\bar C (z_n -b)^{-1} D(z_n)^{-1} Q(z_n)^{-1}$ can be replaced by its 
value at $z_0$ plus an $O(e^{-\veps n})$ error. \eqref{n10.4} implies \eqref{4.12}. 
\end{proof}

\section{Zeros Near Regular Points} \lb{s9}

We call a point $z$ with $\abs{z}=b$ singular if either $z=b$ or $D(1/\bar z)^{-1} 
=0$. Regular points are all not singular points on $\{z\mid\abs{z}=b\}$. There are 
at most a finite number of singular points. In this section, we will analyze zeros 
of $\varphi_n(z)$ near regular points. In the next section, we will analyze the 
neighborhood of singular points. 

We will use Rouch\'e's theorem to reduce zeros of $\varphi_n(z)$ to zeros of 
$(z/b)^n -g(z_0 b)$ ($g$ defined in \eqref{9.1}) for $z-z_0 b$ small. Suppose 
$bz_0$ with $z_0\in\partial\bbD$ is a regular point. Define 
\begin{equation} \lb{9.1} 
g(z) = \f{\bar C\, \ol{D(1/\bar z)}}{D(z) (b-z)} 
\end{equation} 
which is regular and nonvanishing at $z_0$, so 
\begin{equation} \lb{9.2} 
g(bz_0) = a e^{i\psi} 
\end{equation} 
with $a>0$ and $\psi\in [0,2\pi)$. Pick $\delta <b\Delta_2^{-1}$ so that 
\begin{equation} \lb{9.3} 
\abs{z-bz_0} <\delta \Rightarrow \abs{g(z) -g(bz_0)} < \f{a}{4} 
\end{equation} 
Define 
\begin{align} 
h_1^{(n)}(z) &= \biggl( \f{z}{b}\biggr)^n - g(z_0 b) \lb{9.4} \\
h_2^{(n)}(z) &= \biggl( \f{z}{b}\biggr)^n - g(z) \lb{9.5} \\
h_3^{(n)} (z) &= \f{\varphi_n (z)\, \ol{D(1/\bar z)}\,}{b^n} \notag 
\end{align} 
Theorem~\ref{T8.1} implies that 
\begin{equation} \lb{9.6} 
\abs{z-bz_0}<\delta \Rightarrow \abs{h_2^{(n)}(z) - h_3^{(n)} (z)} \leq D_2 
\Delta_1^n 
\end{equation} 

\begin{theorem}\lb{T9.1} Let $z_0$ be a regular point and $\delta < b\Delta_2^{-1}$ 
so that \eqref{9.3} holds. Let $j_1 \leq j_2$ be integers with $\abs{j_k} < (n-1) 
\delta/2$. Let 
\begin{equation} \lb{9.7}
I_n =\biggl\{z\biggm| \biggl| \f{\abs{z}}{b} -1\biggr| < \f{\delta}{2b}, \, 
\arg \biggl( \f{z}{z_0}\biggr) \in \biggl( \f{\psi}{n} + \f{2\pi j_1}{n} - 
\f{\pi}{n}, \f{\psi}{n} + \f{2\pi j_2}{n} + \f{\pi}{n} \biggr) \biggr\}
\end{equation} 
Then for $n$ large, $I_n$ has exactly $(j_2 - j_1) +1$ zeros 
$\{z_\ell^{(n)}\}_{\ell=1}^{j_2 - j_1 +1}$ of $\varphi_n$. Moreover, 
\begin{SL} 
\item[{\rm{(a)}}] 
\begin{align} 
\abs{z_\ell^{(n)}} &= b\biggl( 1 + \f{1}{n} \, \log \abs{g(z_\ell^{(n)})} 
+O \biggl( \f{1}{n^2}\biggr) \biggr) \lb{9.8} \\
&= b \biggl( 1 + \f{1}{n} \, \log a + O\biggl( \f{\delta}{n}\biggr) + 
O \biggl( \f{1}{n^2}\biggr)\biggr) \lb{9.9} 
\end{align} 
\item[{\rm{(b)}}] 
\begin{align} 
\arg z_\ell^{(n)} &= \f{\arg g(z_k^{(n)})}{n} + \f{2\pi \ell}{n} + 
O\biggl( \f{1}{n^2}\biggr) \lb{9.10} \\
&= \f{\psi}{n} + \f{2\pi \ell}{n} + O\biggl( \f{\delta}{n}\biggr) + 
O\biggl( \f{1}{n^2}\biggr) \lb{9.11} 
\end{align} 
\item[{\rm{(c)}}] If $\abs{j_k}\leq J$\!, \eqref{9.9} and \eqref{9.11} hold 
without the $O(\delta/n)$ term.  
\end{SL} 
\end{theorem} 

\begin{proof} Note first that in $I_n$, since $\abs{\abs{z}-b} \leq \f{\delta}{2}$ 
and $\abs{\arg(z/z_0)} \leq \f{\delta}{2}$, we have $\abs{z-z_0 b} < \delta$. Consider 
the boundary of the region $I_n$ which has two arcs at $\abs{z}=b (1\pm \delta)$ 
and two straight edges are $\arg (\f{z}{z_0}) = \f{\psi}{n} + \f{2\pi j_1}{n} - 
\f{\pi}{n}$ and $\arg (\f{z}{z_0}) = \f{\psi}{n} + \f{2\pi j_2}{n} + \f{\pi}{n}$. 

We claim on $\partial I_n$, we have for $n$ large that 
\begin{equation} \lb{9.12} 
\abs{h_j(z) - h_1(z)} \leq \tfrac12\, \abs{h_1(z)} \qquad j=2,3 
\end{equation} 
Consider the $4$ pieces of $\partial I_n$: 

\smallskip 
\noindent{\ul{$\abs{z}=b(1+\delta)$}.} \quad $\abs{h_1^{(n)} (z)} \geq (1+\delta)^n 
-a > a/2$ for $n$ large so, by \eqref{9.3} in this region, 
\[
\abs{h_1^{(n)} - h_2^{(n)}} \leq \tfrac{a}{4} < \tfrac12\, \abs{h_1^{(n)}} 
\]
and clearly, by \eqref{9.6} for $n$ large, 
\[
\abs{h_1^{(n)} - h_3^{(n)}} < \tfrac12\, \abs{h_1^{(n)}} 
\]

\smallskip 
\noindent{\ul{$\abs{z}=b(1-\delta)$}.} \quad $\abs{h_1^{(n)}(z)} \geq a - (1-\delta)^n 
> a/2$ for $n$ large so, as above, \eqref{9.12} holds there. 

\smallskip 
\noindent{\ul{$\abs{\arg z - \f{\psi}{n} - \f{2\pi j_\ell}{n}} = \f{\pi}{n}$}.} \quad 
$z^n = -e^{i\psi} \abs{z}^n$ so 
\[ 
\abs{h_1^{(n)}(z)} = \biggl| \f{z}{b}\biggr|^n + a >a 
\]
and the argument follows the ones above to get \eqref{9.12}. Thus, \eqref{9.12} holds.  

\smallskip 
By Rouch\'e's theorem, $h_1^{(n)}, h_2^{(n)}, h_3^{(n)}$ have the same number of zeros 
in each $I_n$. By applying this to each region with $j_1 =j_2$ and noting that 
$h_1^{(n)}$ has exactly one zero in such a region, we see that each $h_k^{(n)}$ has 
$j_1 + j_2 + 1$ zeros in each $I_n$, and there is one each in each pie slice of 
angle $2\pi/n$ about angles $\psi/n + 2\pi\ell/n$. 

At the zeros of $h_3^{(n)}$, we have $\abs{h_2^{(n)}(z)}\leq D_2 \Delta_1^n$ so 
\begin{align*} 
\biggl| \f{z_\ell^{(n)}}{b}\biggr| &= e^{\log \abs{g(z_\ell^{(n)})}/n} + 
O(\Delta_2^n) \\ 
&= 1 + \f{\log \abs{g(z_\ell^{(n)})}}{n} + O\biggl( \f{1}{n^2}\biggr) 
\end{align*} 
proving \eqref{9.8}. The proof of \eqref{9.10} is similar. Since $\abs{g(z) - 
g(z_0)} \leq C\delta$, we get \eqref{9.9} and \eqref{9.11}. If $\abs{j_k}\leq J$, 
$\abs{g(z_\ell^{(n)}) - g(z_0)} \leq O(1/n)$ so the $O(\delta/n)$ term is not 
needed, proving (c). 
\end{proof} 

\begin{theorem}\lb{T9.2} In each region $I_n$ of Theorem~\ref{T9.1}, the zeros 
$z_\ell^{(n)}$ obey 
\begin{align} 
\abs{z_\ell^{(n)}} - \abs{z_{\ell+1}^{(n)}} &= O\biggl( \f{1}{n^2}\biggr) \lb{9.14} \\
\biggl| \, \arg z_{\ell+1}^{(n)} - \arg z_\ell^{(n)} - \f{2\pi}{n}\biggr| 
&= O\biggl( \f{1}{n^2}\biggr) \lb{9.15} 
\end{align} 
\end{theorem} 

\begin{proof} By \eqref{9.9} and \eqref{9.10}, we have 
\[
\abs{z_{\ell+1}^{(n)} - z_\ell^{(n)}} \leq \f{C}{n} 
\]
(in fact, $C=2\pi b+ O(\delta)$). Thus 
\[
\abs{g(z_{\ell+1}^{(n)}) - g(z_\ell^{(n)})} \leq \f{C_1}{n} 
\]
so \eqref{9.8} and \eqref{9.10} imply \eqref{9.14} and \eqref{9.15}. 
\end{proof}

\section{Zeros Near Singular Points} \lb{s10} 

We first consider the singular point $z=b$ which is always present, and then 
turn to other singular points which are quite different: 

\begin{theorem} \lb{T10.1} Let $\{\alpha_j\}_{j=0}^\infty$ be a sequence of Verblunsky 
coefficients obeying BLS asymptotics \eqref{1.9}. Fix any positive integer, $j$, with 
$j <(n-1) \delta/2$ where $\delta < b\Delta_2^{-1}$ is picked so that 
\begin{equation} \lb{10.1} 
\abs{z-b} < \delta \Rightarrow \abs{g(z) -1} < \tfrac14 
\end{equation} 
Then 
\[
I_n = \biggl\{ z\biggm| \biggl| \f{\abs{z}}{b} -1 \biggr| < \f{\delta}{2b}, \, 
\abs{\arg (z)} < \f{2\pi j}{n} + \f{\pi}{n} \biggr\} 
\]
has exactly $2j$ zeros $\{z_\ell^{(n)}\}_{\ell=1}^j \cup \{z_\ell^{(n)}\}_{\ell= -j}^{-1}$ 
with 
\begin{align} 
\abs{\abs{z_\ell^{(n)}}-b} &= O\biggl( \f{\delta}{n}\biggr) + 
O\biggl( \f{1}{n^2}\biggr) \lb{10.2} \\ 
\arg z_\ell^{(n)} &= \f{2\pi \ell}{n} + O\biggl( \f{\delta}{n}\biggr) + 
O\biggl( \f{1}{n^2}\biggr) \lb{10.3} 
\end{align} 
Moreover, for each fixed $\ell=\pm 1, \pm 2, \dots$, 
\begin{equation} \lb{10.4} 
\abs{z_\ell^{(n)} -be^{2\pi i\ell/n}} = O\biggl( \f{1}{n^2}\biggr)
\end{equation} 
\end{theorem} 

{\it Remark.} We emphasize again the zero at $z=b$ is ``missing," that is, 
$z_\ell^{(n)}$, $\abs{\ell}\geq 2$, has its nearest zeros in each direction a 
distance $b 2\pi/n + O(1/n^2)$, while $z_\pm^{(n)}$ has a zero on one side at 
this distance, but on the other at distance 
$b 4\pi/n + O(1/n^2)$. 

\begin{proof} This is just the same as Theorem~\ref{T9.1}! $g(z)$ is regular at 
$z=b$. Indeed, $g(b)=1$ since $\ol{D(1/\bar z)}^{-1}$ and $\bar CD(z)/(z-b)$ have 
to have precisely cancelling poles. $D(1/\bar z)$ vanishes at $z=b$, so $\varphi_n(z) 
\ol{D(1/\bar z)}/b^n$ which, by the argument of Theorem~\ref{T9.1}, has $2j+1$ zeros 
in $I_n$, has one at $z=b$ from $D(1/\bar z)$ and $2j$ from $\varphi_n(z)$, 
\eqref{10.2}--\eqref{10.4} follow as did \eqref{9.9}, \eqref{9.11}, and 
\eqref{9.14}/\eqref{9.15}. 
\end{proof} 

By a compactness argument and Theorems~\ref{T9.1} and \ref{T10.1}, we have the following 
strong form of Theorem~\ref{T4.3} when there are no singular points other than $z=b$:  

\begin{theorem}\lb{T10.2} Let $\{\alpha_j\}_{j=0}^\infty$ be a sequence of Verblunsky 
coefficients obeying the BLS condition \eqref{1.9}. Suppose $z=b$ is the only singular 
point, that is, $D(z)^{-1}$ is nonvanishing on $\abs{z}=b^{-1}$. Let $\{w_j\}_{j=1}^J$ 
be the Nevai-Totik zeros and $m_j$ their multiplicities. Let $\delta$ be such that for 
all $j\neq k$, $\abs{w_j - w_k} >2\delta$ and $\abs{w_j} > b+2\delta$, and let $M= 
\sum_{j=1}^J m_j$ the total multiplicity of the NT zeros. Let $C> \sup_{\abs{z}=b} 
\abs{\log(\abs{g(z)})}$. Then 
\begin{SL} 
\item[{\rm{(1)}}] For some large $N$ and all $n>N$\!, the zeros of $\varphi_n(z)$ are 
$m_j$ zeros in $\{z\mid\abs{z-w_j}<\delta\}$ for $j=1, \dots, J$ and $n-M$ in the annulus 
\begin{equation} \lb{10.5} 
\abs{\,\abs{z}-\abs{b}\,}  < \f{C}{n}
\end{equation} 

\item[{\rm{(2)}}] The $n-M$ zeros in the annulus can be labelled by increasing arguments 
$0<\arg (z_1^{(n)}) <\cdots < \arg(z_{n-J}^{(n)}) <2\pi$, and with $\arg(z_{n-J+1}^{(n)}) 
= \arg (z_1^{(n)}) + 2\pi$ and $\abs{z_{n-J+1}^{(n)}} = \abs{z_1^{(n)}}$, we have 
\begin{equation} \lb{10.6} 
\arg (z_{k+1}^{(n)}) - \arg (z_k^{(n)}) = \f{2\pi}{n} + O\biggl( \f{1}{n^2}\biggr)  
\end{equation} 
and 
\begin{equation} \lb{10.7} 
\f{\abs{z_{k+1}^{(n)}}}{\abs{z_k^{(n)}}} = 1 + O\biggl( \f{1}{n^2}\biggr) 
\end{equation} 
uniformly in $k=1,2, \dots, n-J+1$.  
\end{SL} 
\end{theorem} 

{\it Remarks.} 1. By \cite{MGH}, if $\alpha_n =-b^n$, all zeros have $\abs{z}=b$, 
so it can happen that there are no $O(1/n)$ terms in $\abs{z_k^{(n)}}$. However, 
since there are $n$ zeros in $(\f{2\pi}{n}, \f{n-1}{n} 2\pi)$ (i.e., the zero near 
$z=b$ is missing), there are always either NT zeros or $O(1/n^2)$ corrections in 
some $\arg (z_n^{(\ell)})$. 

\smallskip 
2. While it can happen that there is no $O(1/n)$ term, its absence implies strong 
restrictions on $\alpha_n$. For 
\begin{equation} \lb{10.9a} 
\abs{g(z)}=1 \quad\text{on}\quad \abs{z}=b 
\end{equation} 
implies $\ol{g(b/\bar z)}\, g(z) =1$ near $\abs{z}=b$, and that equation allows 
analytic continuations of $g$, and so of $D$ or $D^{-1}$. In fact, \eqref{10.9a} 
implies that $D(z)^{-1}$ is analytic in $\{z\mid \abs{z}<b^{-3}\}$ except for 
a pole at $z=b^{-1}$ and that implies, by Theorem~7.2.1 of \cite{OPUC1}, that 
$\alpha_n =-Cb^n + O(b^{2n})(1-\veps)^n$ for all $\veps$. Thus, if the BLS 
condition holds but $\liminf \abs{\alpha_n -Cb^n}^{1/n} >b^2$, then there 
must be $O(1/n)$ corrections to $\abs{z_k^{(n)}}=b$. Note that for the 
Roger's Szeg\H{o} polynomials, the poles of $D(z)^{-1}$ are precisely at 
$z\in\{b^{-2k-1}\}_{k=0}^\infty$ (see \cite[Eqn.~(1.6.59)]{OPUC1}), 
consistent with the $\abs{z}<b^{-3}$ statement above.  

\smallskip 
We now turn to an analysis of the other singular points. As a warmup, 
we study zeros of 
\begin{equation} \lb{10.8} 
f_n(z) =z^n -K(1-z)^k 
\end{equation} 
where 
\[ 
K=ae^{i\psi} 
\]
is nonzero, $a>0$, $k$ a fixed positive integer, and we take $n\to\infty$.  

We begin by localizing $\abs{z}$ and $\abs{z-1}$. 

\begin{proposition}\lb{P10.3} 
\begin{SL} 
\item[{\rm{(i)}}] There are $M>0$ and $N_0$ so that if $\abs{z} \geq 1 + M/n$ and 
$n\geq N_0$, then $f_n(z) \neq 0$. 
\item[{\rm{(ii)}}] There is $N_1$ so if $n\geq N_1$ and 
\begin{equation} \lb{10.9} 
\abs{z} \leq 1 -2k \,\f{\log n}{n} 
\end{equation}
then $f_n(z) \neq 0$. 
\item[{\rm{(iii)}}] There is $N_2$ so if $n\geq N_2$ and 
\begin{equation} \lb{10.10} 
\abs{z-1} \leq \f{k}{2} \, \f{\log n}{n} 
\end{equation}
then $f_n(z) \neq 0$. 
\end{SL}  
\end{proposition} 

{\it Remark.} If one proceeds formally and lets $y=1-z$, then $f_n(z) =0$ is 
equivalent to $(1-y)^n =Ky^k$ and finds  
\begin{align*} 
y &= \f{k}{n}\, (-\log y) + O(y^2) + \f{\log a}{n} + \f{i\psi}{n} \\
&= \f{k}{n} \, \log n - \f{k}{n} \, \log k - \f{k}{n} \, \log (\log n) + 
O\biggl(\biggl( \f{\log n}{n}\biggr)^2\biggr) 
\end{align*} 
It was this formal calculation that caused us to pick $2k \f{\log n}{n}$ in 
\eqref{10.9} and $\f{k}{2} \f{\log n}{n}$ in \eqref{10.10}. 

\begin{proof} (i) For $n$ sufficiently large and $\abs{z}\geq 2$, $\abs{z}^n \geq 
a(1+\abs{z})^k$ and for $n$ large, $( 1+\f{M}{n})^n \geq e^{M/2} > a2^k$ if $M$ 
is suitable, so for such $M$ and $n$, $\abs{z}^n \geq a(1+\abs{z})^k$ for 
$1+\f{M}{n} \leq\abs{z}\leq 2$. 

\smallskip 
(ii) If \eqref{10.9} holds, 
\[
\abs{1-z}^k \geq (1-\abs{z})^k \geq (2k)^k \biggl( \f{\log n}{n}\biggr)^k 
\]
while $\log (1+x) \leq x$ implies 
\begin{align*} 
\abs{z}^n &= \exp (n\log (1+ (\abs{z}-1))) \\
&\leq \exp (-n(1-\abs{z})) \\
&\leq \exp (-2k \log n) =n^{-2k} 
\end{align*} 
so for $n$ large, $a\abs{1-z}^k >\abs{z}^n$. 

\smallskip 
(iii) If \eqref{10.10} holds, then  
\[
\abs{z} \geq 1-\abs{z-1} \geq 1 - \f{k}{2} \, \f{\log n}{n} 
\]
so for $n$ large, 
\begin{align*} 
\abs{z}^n &\geq \exp\biggl(\f{3}{4}\, k \log n \biggr)\\ 
&= n^{-\f{3}{4} k} \\
&> \biggl( \f{k}{2} \, \f{\log n}{n}\biggr)^k \\
&\geq \abs{1-z}^k 
\end{align*} 
so $f_n(z)\neq 0$. 
\end{proof} 

With this, we are able to control ratios of lengths of nearby zeros.  

\begin{proposition}\lb{P10.4} 
\begin{SL} 
\item[{\rm{(i)}}] If $z,w$ are two zeros of $f_n(z)$ and 
\begin{equation} \lb{10.11} 
\biggl| \, \arg \biggl( \f{z}{w}\biggr)\biggr| \leq \f{C_0 \log n}{n}  
\end{equation} 
then, for $n$ large, 
\begin{equation} \lb{10.12} 
\biggl| \,\biggl| \f{z}{w}\biggr| -1 \, \biggr| \leq \f{D_0}{n}  
\end{equation} 

\item[{\rm{(ii)}}] If $z,w$ are two zeros of $f_n(z)$ and 
\begin{equation} \lb{10.13} 
\biggl|\, \arg \biggl( \f{z}{w}\biggr)\biggr| \leq \f{C_1}{n}   
\end{equation} 
then 
\begin{equation} \lb{10.14} 
\biggl| \,\biggl| \f{z}{w}\biggr| -1\, \biggr| \leq \f{D_1}{n\log n}   
\end{equation}
\end{SL} 
\end{proposition} 

{\it Remark.} In both cases, $D_j$ is a function of $C_j$, $K$\!, and $k$. 

\begin{proof} Since $z$ and $w$ are both zeros, 
\[ 
\biggl(\f{z}{w}\biggr)^n = \f{(z-1)^k}{(w-1)^k} 
\]
so, by (iii) of the last proposition, $\abs{w-1}\geq \f{k}{2}\, \f{\log n}{n}$. Thus  
\begin{align*} 
\biggl| \f{z}{w}\biggr| &\leq \biggl( 1 + \f{\abs{z-w}}{\abs{w-1}}\biggr)^{k/n} \\
&\leq \biggl( 1+ \f{2n\abs{z-w}}{\log n} \biggr)^{k/n} 
\end{align*} 
If \eqref{10.11} holds, $\abs{\f{z}{w}} \leq (1+\f{2C_0}{n})^{k/n} \leq 1 + 
\f{\ti D_0}{n}$ for $n$ large. Interchanging $z$ and $w$ yields \eqref{10.12}. 
The argument from \eqref{10.13} to \eqref{10.14} is identical. 
\end{proof} 

As the final step in studying zeros of \eqref{10.8}, we analyze arguments of nearby zeros. 

\begin{proposition} \lb{P10.5} Let $z_0$ be a zero of $f_n(z)$. Fix $C_2$. Then 
\begin{SL} 
\item[{\rm{(i)}}] If $w$ is also a zero and $\abs{w-z_0} \leq C_2/n$, then for 
some $\ell\in\bbZ$, 
\begin{equation} \lb{10.15} 
\arg\biggl( \f{w}{z}\biggr) = \f{2\pi\ell}{n} + O\biggl( \f{1}{n\log n}\biggr) 
\end{equation}
with the size of the error controlled by a $C_2$, $K$\!, and $k$. 

\item[{\rm{(ii)}}] For each $L$ and $n$ large, there exists exactly one zero 
obeying \eqref{10.15} for each $\ell=0,\pm 1, \pm2, \dots, \pm L$. 

\item[{\rm{(iii)}}] For any $\psi_0\in [0,2\pi)$ and $\delta$, there is $N$ so 
for $n\geq N$\!, there is a zero of $f_n$ with $\arg z\in (\psi_0 -\f{\pi+\delta}{n}, 
\psi_0 + \f{\pi +\delta}{n})$. 
\end{SL} 
\end{proposition} 

{\it Remarks.} 1. These propositions imply that the two nearest zeros to $z_0$ are  
$z_0 e^{\pm 2\pi i/n} + O(1/n\log n)$. 

\smallskip 
2. If $k$ is odd, the argument of $(z-1)^k$ changes by $\pi k$ as $z$ moves 
through $1$ along the circle. Thus there are $O(1)$ shifts as zeros swing around 
the forbidden circle $\abs{z-1}\sim \f{k}{n} \log n$. Since there are $\log n$ zeros 
near that circle, the phases really do slip by $O(1/n\log n)$. This is also why we 
do not try to specify exact phases, only relative phases. 

\begin{proof} Since $z_0^n =K(z_0-1)^k$, we have 
\begin{equation} \lb{10.16} 
\f{f(z)}{z_0^n} = \biggl( \f{z}{z_0}\biggr)^n - \f{(z-1)^k}{(z_0-1)^k} 
\end{equation} 
If $\abs{z-z_0}\leq C_2/n$, then since $\abs{z_0-1}\geq \f{1}{2} K \f{\log n}{n}$ 
(by Proposition~\ref{P10.3}(iii)),  
\begin{equation} \lb{10.17} 
\biggl| \f{(z-1)^k}{(z_0-1)^k} -1 \biggr| \leq \f{D\abs{z-z_0}}{\abs{z_0-1}} 
\leq \f{D_2}{\log n} 
\end{equation} 

Fix $\eta$ and any $z_0$ with $\abs{z_0-1}\geq \f12 K \f{\log n}{n}$, we want to look 
at solutions of 
\begin{equation} \lb{10.18} 
\biggl( \f{z}{z_0}\biggr)^n = e^{i\eta} \, \f{(z-1)^k}{(z_0-1)^k}  
\end{equation} 
with $\abs{z-z_0} \leq C_2/n$. If $f(z_0) =0$ and $\eta =0$, we have solutions of 
\eqref{10.18} are exactly zeros of $f$. 

By \eqref{10.17}, \eqref{10.18} implies 
\[
n(\arg(z) - \arg(z_0)) = \eta + O\biggl( \f{1}{\log n}\biggr) 
\]
whose solutions are precisely 
\begin{equation} \lb{10.19} 
\arg\biggl(\f{z}{z_0}\biggr) = \f{2\pi\ell}{n} + \eta + O\biggl( \f{1}{\log n}\biggr) 
\end{equation}
which is \eqref{10.15} when $\eta =0$. If we prove that for $n$ large, there is 
exactly one solution of \eqref{10.19} for $\ell = 0, \pm 1, \pm 2, \dots, \pm L$, 
we have proven (ii). We also have (iii). For given $\psi_0$, let 
\begin{equation} \lb{10.20} 
h(r) = r^n - \abs{K}\, \abs{re^{i\psi_0} -1}^k  
\end{equation} 
If $n>k$, $h(0) <0$ and $h(\infty) >0$, so there is at least one $r_0$ solving 
$h(r) =0$. Let $z_0 =re^{i\psi_0}$ and define $\eta$ by 
\[
z_0^n = Ke^{-i\eta} (1-z_0)^k 
\]
Solutions of $f_n(z)=0$ are then precisely solutions of \eqref{10.18}, and we have 
the existence statement in (iii) if we prove existence and uniqueness of solutions 
of \eqref{10.18}. 

Fix $M$\!. Consider the following contour with four parts: 
\begin{align*} 
C_1 &= \biggl\{z\biggm| \abs{z}=\abs{z_0}\biggl( 1+ \f{M}{n}\biggr);\, 
\biggl( \eta- \f{\pi}{2}\biggr)\f{1}{n} \leq \arg \biggl( \f{z}{z_0}\biggr) 
\leq \biggl( \eta + \f{\pi}{2}\biggr) \f{1}{n}\biggr\} \\
C_2 &= \biggl\{ z\biggm| \arg\biggl( \f{z}{z_0}\biggr) = \biggl( \eta + \f{\pi}{2} 
\biggr) \f{1}{n}; \, \abs{z_0} \biggl( 1-\f{M}{n}\biggr) \leq \abs{z} \leq 
\abs{z_0} \biggl( 1+\f{M}{n}\biggr) \biggr\} \\
C_3 &= \biggl\{z\biggm| z=\abs{z_0}\biggl( 1-\f{M}{n}\biggr);\, 
\biggl( \eta - \f{\pi}{2}\biggr) \f{1}{n} \leq \arg \biggl( \f{z}{z_0}\biggr) 
\leq \biggl( \eta + \f{\pi}{2}\biggr) \f{M}{n}\biggr\} \\
C_4 &= \biggl\{z\biggm| \arg \biggl( \f{z}{z_0}\biggr) = \biggl( \eta - \f{\pi}{2} 
\biggr) \f{1}{n}; \, \abs{z_0} \biggl( 1-\f{M}{n}\biggr) \leq \abs{z} \leq 
\abs{z_0}\biggl( 1+\f{M}{n}\biggr)\biggr\} 
\end{align*} 
($\pi/2$ can be replaced by any angle in $(0,\pi)$). For $n$ large, $q(z)\equiv 
(z/z_0)^n e^{-i\eta}$ follows arbitrarily close to 
\begin{align*} 
q[C_1] &\cong \biggl\{w\biggm| \abs{w}=e^M;\, \arg (w)\in 
\biggl( -\f{\pi}{2}\,, \f{\pi}{2}\biggr)\biggr\} \\
q[C_2] &\cong \biggl\{w\biggm| e^{-M}\leq \abs{w}\leq e^M;\, 
\arg (w)=\f{\pi}{2}\biggr\} \\
q[C_3] &\cong \biggl\{w\biggm| \abs{w} = e^{-M};\, 
\arg(w)\in \biggl( -\f{\pi}{2}\, , \f{\pi}{2}\biggr)\biggr\} \\
q[C_4] &\cong \biggl\{w\biggm| e^{-M} \leq \abs{w}\leq e^M; \, 
\arg(w)= - \f{\pi}{2}\biggr\} 
\end{align*} 
which surrounds $w=1$ once. 

By \eqref{10.17}, if $\ti q(z)=e^{-i\eta} (z/z_0)^n - (z-1)^k/(z_0-1)^k$, 
then $\ti q[C]$ surrounds $w=0$ once. It follows that \eqref{10.18} has one 
solution inside $C$. This proves existence and uniqueness within $C$. By part (i), 
there are no other solutions with $\abs{z-z_0} < C/n$. 
\end{proof} 

These ideas immediately imply 

\begin{theorem} \lb{T10.6} Let $\{\alpha_j\}_{j=0}^\infty$ be a sequence of 
Verblunsky coefficients obeying BLS asymptotics \eqref{1.9}. Let $bz_0$ be a 
singular point so that $\ol{D(1/\bar z)}^{-1}$ has a zero of order $k$ at $bz_0$. 
Then there is $\delta >0$ with $\delta < b\Delta_2^{-1}$ and $N$ large so that 
for $n\geq N_0$, all zeros of $\varphi_n(z)$ in 
\begin{equation} \lb{10.21} 
S\equiv \biggl\{z\bigg| \Delta_2^{-1} < \biggl| \f{z}{b}\biggr| < \Delta_2; \, 
\biggl| \arg \biggl( \f{z}{z_0}\biggr)\biggr| <\delta \biggr\}
\end{equation}
obey 
\begin{equation} \lb{10.22} 
1-\f{D}{n} \leq \biggl| \f{z}{b}\biggr| \leq 1 + 2k \, \f{\log N}{n} 
\end{equation} 
for some $D$. Moreover, for each $z_1\in S$ and $n>N$\!, there is a zero, $z$, 
in $S$ with 
\begin{equation} \lb{10.23} 
\biggl| \arg\biggl( \f{z}{z_1}\biggr)\biggr| \leq \f{2\pi}{n} 
\end{equation} 

If $z_\ell$ is a zero in $S$, then the two nearest zeros to $z_\ell$ obey 
\begin{align} 
\arg\biggl( \f{z_{\ell\pm 1}}{z_\ell}\biggr) &= \pm \f{2\pi}{n} + 
O\biggl( \f{1}{n\log n}\biggr) \lb{10.24} \\
\biggl| \f{z_{\ell+1}}{z_\ell}\biggr| &= 1 + O\biggl( \f{1}{n\log n}\biggr) \lb{10.25} \\
\abs{bz_0 -z_\ell} &\geq \f{k}{2}\, \f{\log n}{n} \lb{10.26} 
\end{align} 
\end{theorem} 

{\it Remark.} $2k$ in \eqref{10.22} can be replaced by any number strictly bigger than 
$k$ (if we take $N_0$ large enough). Similarly in \eqref{10.23}, $2\pi$ can be any 
number strictly bigger than $\pi$ and $k/2$ in \eqref{10.26} can be any number 
strictly less than $k$. 

\begin{proof} This follows by combining the analysis of $f_n$ above with the ideas 
used to prove Theorem~\ref{T9.1}, picking $h_1^{(z)}=(b/z)^n - K(z-z_0 b)^k$, 
$h_2^{(z)} = (b/z)^n -g(z)^{-1}$, and $h_3^{(z)} =(z-b) D(z) \varphi_n(z)/\bar Cz^n$. 
\end{proof} 

By compactness of $b[\partial\bbD]$ and Theorems~\ref{T9.1}, \ref{T10.1}, and 
\ref{T10.6}, we get the following precise form of Theorem~\ref{T4.3}: 

\begin{theorem} \lb{T10.7} Let $\{\alpha_j\}_{j=0}^\infty$ be a sequence of Verblunsky 
coefficients obeying BLS asymptotics \eqref{1.9}. The total multiplicity $J$ of 
Nevai-Totik zeros is finite. There exists $K$\!, $D>0$, so that for all $n$ large, 
the $n-J$ zeros not near the Nevai-Totik points can be labelled $z_j^{(n)} = 
\abs{z_j^{(n)}} e^{i\theta_j^{(n)}}$, $j=1, \dots, n-J$\!, with $0=\theta_0^{(n)} < 
\theta_1^{(n)} < \cdots < \theta_{n-J+1}^{(n)} =2\pi$ so that  
\begin{SL} 
\item[{\rm{(1)}}] 
\[
1-\f{D}{n} \leq \inf_{j=1, \dots, n-J}\, \f{\abs{z_j^{(n)}}}{b} \leq 
\sup_{j=1, \dots, n-J}\, \f{\abs{z_j^{(n)}}}{b} \leq 1 + K\, \f{\log n}{n} 
\]
\item[{\rm{(2)}}] 
\[
\sup_{j=0, \dots, n-J}\, n \biggl| \theta_{j+1}^{(n)} - \theta_j^{(n)} - 
\f{2\pi}{n}\biggr| = O\biggl( \f{1}{\log n}\biggr) \\
\] 
\item[{\rm{(3)}}] Uniformly in $j$ in $0, \dots, n-J$\!, 
\[ 
\f{\abs{z_{j+1}^{(n)}}}{\abs{z_j}} = 1 + O\biggl( \f{1}{\log n}\biggr) 
\]
\end{SL}
where $\abs{z_0^{(n)}}$ and $\abs{z_{n-J+1}^{(n)}}$ are symbols for $b$. Moreover, 
for $L$ fixed, uniformly in $\ell =1, \dots, L, n-J, n-J-1, \dots, n-J-L+1$, 
\[
z_\ell^{(n)} = \begin{cases} be^{2\pi i\ell/n} & \ell=1, \dots, L \\
be^{-2\pi i(n-J-\ell+1)} & \ell =n-J, n-J-1, \dots, n-J-L+1 
\end{cases} 
\]
\end{theorem}

\section{Comments} \lb{s11} 

One question the reader might have is whether singular points other than $z=b$ ever 
occur and whether there might be a bound on $k$. In fact, there are no restrictions 
for 

\begin{proposition} \lb{P13.1} If $f(z)$ is nonvanishing and analytic in a neighborhood 
of $\bar\bbD$ and 
\begin{equation} \lb{13.1} 
\int \abs{f(e^{i\theta})}^2 \, \f{d\theta}{2\pi} =1
\end{equation} 
then there is a measure $d\mu$ in the Szeg\H{o} class with $D(z;d\mu) =f(z)$. 
\end{proposition} 

\begin{proof} Just take $d\mu =\abs{f(e^{i\theta})}^2 \f{d\theta}{2\pi}$. 
\end{proof} 

\begin{example}\lb{E13.2} Pick disjoint points $z_1, \dots, z_\ell$ with $\arg(z_j)\neq 0$ 
and $\abs{z_j}=b^{-1}$ and positive integers $k_1, \dots, k_\ell$. Let 
\[
f(z) = c \biggl[\, \prod_{j=1}^\ell (z-z_j)^{-k_j}\biggr] (z-b^{-1}) 
\]
where $c$ is chosen so that \eqref{13.1} holds. Then $D^{-1}$ has a single pole at 
$b^{-1}$ and so, by Theorem~\ref{T1.2}, $\alpha_j (d\mu)$ obeys the BLS condition. 
By Proposition~\ref{P13.1}, 
\[
D(z)^{-1} = c^{-1} \, \f{1}{z-b^{-1}} \, \prod_{j=1}^\ell (z-z_j)^{k_j} 
\]
and so has zeros at $z_j$ of order $k_j$. Thus we have singular points at 
$1/\bar z_j$ of order $k_j$.  

Let me emphasize that this procedure also lets one move Nevai-Totik zeros without changing 
the leading asymptotics of the Verblunsky coefficients or the error estimate. 
\qed 
\end{example} 

Secondly, we want to prove that in a suitable generic sense, the only singular point 
is $z=b$. Fix $b$ and $\Delta$ in $(0,1)$. Let $\calV_{b,\Delta}$ be the space of 
sequences $\alpha_j$ in $\bigtimes\bbD$ so that for some $C\neq 0$, 
\begin{equation} \lb{13.2} 
\abs{C} + \sum\, \abs{(\alpha_j -Cb^j)\Delta^{-j}b^{-1}} \equiv \|\alpha\| < \infty 
\end{equation} 
Obviously, any $\alpha\in\calV_{b,\Delta}$ obeys the BLS condition. 

\begin{theorem}\lb{T13.3} $\{\alpha\in\calV_{b,\Delta}\mid\text{The only single 
point is $z=b$}\}$ is a dense open set of $\calV_{b,\Delta}$. 
\end{theorem} 

\begin{proof} The proof of Theorem~7.2.1 of \cite{OPUC1} shows that the finiteness 
condition of \eqref{13.2} is equivalent to $D(z)^{-1} -C(z-b^{-1})^{-1}$ lying 
in the Wiener class of the disk $\{z\mid\abs{z}\leq (\Delta b)^{-1}\}$ and that 
$\|\cdot\|$ convergence is equivalent to Wiener convergence of $(z-b^{-1}) 
D(z)^{-1}$. 

That the set of no singular points other than $z=b$ is open follows from the fact 
that the just mentioned Wiener convergence implies uniform convergence, and so 
$D(z)^{-1}$ is nonvanishing on $\{z\mid \abs{z}= b^{-1}\}$ is an open set. 

As noted in Example~\ref{E13.2}, it is easy to move zeros by multiplying 
$D(z)^{-1}$ by $\f{z-z_1}{z-z_2}$. This shows that any $D(z)^{-1}$ in the Wiener 
space is a limit of such $D$'s nonvanishing on $\{z\mid\abs{z}=b^{-1}\}$ and shows 
that the set of nonsingular $\alpha$ is dense. 
\end{proof} 

It is a worthwhile and straightforward exercise to compute the change of $\arg(g)$ 
along $\{z\mid \abs{z}=b\}$ and so verify that the number of zeros on the critical 
circle is exactly $n$ minus the total order of the Nevai-Totik zeros.

\bigskip


\end{document}